\documentclass[a4paper,10pt]{article}
\usepackage{amsfonts}
\newcommand{\field}[1]{\mathbb{#1}}

\title{Conformally parametrized surfaces associated with $\field{C}P^{N-1}$ 
sigma models}

\author{\vspace{1cm}\\
         {\bf A.\ M.\ Grundland}$^{1,2}$
          \thanks{E-mail address:
           grundlan@crm.umontreal.ca} \,, 
            {\bf W.\ A.\ Hereman}$^3$
             \thanks{E-mail address:
                  whereman@mines.edu} 
            {\,, and \bf \.{I}.\ Yurdu\c{s}en}$^1$
             \thanks{E-mail address:
       yurdusen@crm.umontreal.ca}
          \\
          \\$^1$Centre de Recherches Math\'{e}matiques, Universit\'{e} 
                   de Montr\'{e}al, 
               \\ CP 6128, Succ.\ Centre-Ville, Montr\'{e}al, 
                        Qu\'{e}bec H3C 3J7, Canada
                   \\
                 \\$^2$Universit\'{e} du Qu\'{e}bec, 
                     Trois-Rivi\`{e}res, CP500, QC, G9A 5H7, Canada
                    \\
                  \\$^3$Department of Mathematical and Computer Sciences,
                     \\ Colorada School of Mines, 
                   Golden, CO, 80401-1887, U.S.A.}

\date{\today}
\begin{document}

\maketitle

\begin{abstract}
Two-dimensional conformally parametrized surfaces immersed in 
the $su(N)$ algebra are investigated. The focus is on surfaces 
parametrized 
by solutions of the equations for the $\field{C}P^{N-1}$ sigma model. 
The Lie-point symmetries of the $\field{C}P^{N-1}$ model are computed for
arbitrary $N$. The Weierstrass formula for immersion is determined 
and an explicit formula for a moving frame on a surface is constructed. 
This allows us to determine the structural equations and geometrical 
properties of surfaces in $\field{R}^{N^2-1}$. 
The fundamental forms, Gaussian and mean curvatures, Willmore functional 
and topological charge of surfaces are given explicitly in terms of 
any holomorphic solution of the $\field{C}P^{2}$ model. 
The approach is illustrated through several examples, including surfaces 
immersed in low-dimensional $su(N)$ algebras.
\end{abstract}




\section{Introduction \label{intro}}
Group theoretical methods have proven to be very useful for studying 
surfaces immersed in multi-dimensional spaces and for computing their 
main geometric characteristics 
\cite{Bobenko1, Fokas1, Fokas2, Helein1, Dorfmeister}. 
It was shown in 
\cite{Grundlandetal, GrundlandSnobl, Konopelchenko1, Konopelchenko2} 
that the problem of Weierstrass immersion of two-dimensional smooth 
surfaces in multi-dimensional Euclidean spaces is related to the 
surfaces in Lie algebras associated with the $\field{C}P^{N-1}$ models. 
The main feature of this approach is that it allows one to replace the 
methods based on Dirac-type equations by a formalism connected with 
completely integrable $\field{C}P^{N-1}$ models. 
The task of finding an increasing number of surfaces is related to 
choosing a suitable Lie representation of the $\field{C}P^{N-1}$ model. 
Group analysis makes it possible to construct 
algorithms proceeding directly from the equations of the 
$\field{C}P^{N-1}$ model 
and without referring to any additional considerations. 
The techniques for constructing two-dimensional surfaces immersed in 
$su(N)$ algebras, obtained from integrable models, 
are better understood for low-dimensional $\field{C}P^{N-1}$ models. 
In that case, the geometric features of surfaces so obtained are interesting 
and the subject of ongoing study. 
A review of recent developments related to integrable models can be found in 
\cite{Helein2, Guest, Kobayashi, Nomizu2}.

Over the last century and a half, the Weierstrass formula for immersion 
of surfaces in Lie groups, Lie algebras and homogeneous spaces has been 
used extensively in various areas of mathematics, physics, chemistry 
and biology. 
We now list some of the most important examples. 

In mathematics, the topic is of central importance in the formulation 
of the classical theory of surfaces. 
In particular, immersions are useful for studying surfaces with techniques 
of completely integrable continuous and discrete systems, 
as well as for the development and application of numerical tools 
\cite{Bobenko2,Rogers}. 
A description of the monodromy of solutions of Painlev\'{e} equations 
is yet another important application \cite{Bobenko3}.

In physics, the concept has numerous applications in, e.g., 
two-dimensional gravity \cite{Gross}, field and string theory 
\cite{Polchinski, Nelson}, 
statistical physics (e.g., growth of crystals, surface waves, 
dynamics of vortex sheets, the two-body correlation function of the 
two-dimensional Ising model \cite{David}), 
fluid dynamics (e.g., motion of boundaries between regions of 
differing densities and velocities \cite{Sommerfeld}), 
plasma physics (geometry of magnetic surfaces and constant pressure 
surfaces in various fusion devices like tokomaks, stellarators, 
magnetic mirrors \cite{Chen}).

In chemistry, descriptions of energy and momentum transport along a 
polymer molecule constitute a significant area of application for the 
theory of immersions \cite{Davidov, Ou}.
In biology, the theory is frequently used in the study of the model for 
the Canham-Helfrich membrane and its continuous deformations 
\cite{Safram, Landolfi}.

In general, the algebraic approach to the equations describing surface 
immersion has been proven to be very fruitful from a computational point of view. 
In addition, the geometric approach is of primary importance to the 
derivation and characterization of the governing equations for related 
phenomena in physics and other applied sciences.

This paper follows-up on research in \cite{Grundlandetal}, where
surfaces immersed in $su(N+1)$ algebras obtained via $\field{C}P^N$ 
models were investigated. 
We generalize the results and also correct some formulae. 
To be precise, the new results presented in this paper include the 
Lie-point symmetry algebra of the $\field{C}P^{N-1}$ model 
for arbitrary $N$. 
We also give new examples of surfaces immersed in the $su(N)$ algebra 
invariant under the scaling  symmetries whose Gaussian curvature 
always vanishes. 
We delve deeply into the geometrical aspects of surfaces in $su(3)$ 
obtained from the $\field{C}P^2$ model. 
For that case, we identify the moving frame and the structural equations, 
as well as the Willmore functional and the topological charge. 
The main goal of this paper is to provide a comprehensive, self-contained 
approach to the subject.

The paper is organized as follows. 
In Section \ref{preliminary}, we briefly review some basic notions and 
properties concerning the Euler-Lagrange equations associated with 
the $\field{C}P^{N-1}$ models. 
In Section \ref{secweierstrass}, we discuss the Weierstrass formula for 
immersion in connection with the $\field{C}P^{N-1}$ model, derive the 
induced metric and compute the scalar curvature. 
Section \ref{secliepointsymm} is devoted to the Lie-point symmetries 
of the equations of the $\field{C}P^{N-1}$ model for arbitrary $N$. 
Section \ref{secimmsu3cp2} covers the analysis of the immersion of 
surfaces in the $su(3)$ algebra arising from the $\field{C}P^2$ model. 
In Section \ref{Weierstrassformulaforsu1su2} we investigate the Weierstrass 
aspects for immersion of surfaces in the $su(2)$ and $su(3)$ algebras which 
are associated with the $\field{C}P^1$ and $\field{C}P^2$ models, 
respectively.
Section \ref{secexamples} deals with applications of the Weierstrass 
formula for the immersion of surfaces in the $su(2)$ and $su(3)$ algebras, 
as well as surfaces immersed in the $su(N)$ algebra invariant under the 
scaling symmetries.

\section{The Euler-Lagrange equations associated with the 
$\field{C}P^{N-1}$ sigma models \label{preliminary}}
To keep the paper self-contained, 
we briefly review basic notions and properties of the $\field{C}P^{N-1}$ 
sigma models (see {e.g.,} \cite{Helein2, Zakrzewski, Mikhailov} and 
references therein).
The domain of definition for the sigma model is assumed to be an open, 
connected and simply connected set $\Omega \subset \field{C}$ 
with the Euclidean metric
\begin{equation}
ds^2=d\xi d\bar{\xi}=(d\xi^1)^2+(d\xi^2)^2\,, \qquad 
\xi=\xi^1+i \xi^2\,, \label{euclideanmetric}
\end{equation}
where $\xi$ and $\bar{\xi}$ are local coordinates in $\Omega$. 
In the case of the $\field{C}P^{N-1}$ models the target space is 
a $(N-1)$-dimensional complex projective space $\field{C}P^{N-1}$, 
which is defined as the set of all complex lines in $\field{C}^N$. 
The manifold structure on it is defined by an open covering
\begin{equation}
\mathcal{U}_k=\{[z]\,|\,\,z \!\in \field{C}^N, z_k\ne 0\}\,, 
\quad k=1,\ldots,N\,, 
\label{manifoldstr}
\end{equation}
where $[z]=\rm{span}\{z\}$ and the coordinate maps 
$h_k:\mathcal{U}_k\rightarrow \field{C}^{N-1}$ are defined by 
\begin{equation}
h_k(z)=\left(\frac{z_1}{z_k}, \ldots ,\frac{z_{k-1}}{z_k}, 
\frac{z_{k+1}}{z_k}, \ldots , \frac{z_N}{z_k}\right).
\end{equation}
We are interested in maps of the form 
$[z]:\Omega \rightarrow \field{C}P^{N-1}$, which are stationary 
points of the action functional 
\begin{equation}
S=\frac{1}{4}\int_{\Omega}(D_{\mu} z)^{\dagger} (D^{\mu} z) 
d\xi d\bar{\xi}\,, \qquad z^{\dagger}\cdot z=1\,. 
\label{action}
\end{equation}
Here, $D_{\mu}$ and $D^{\mu}$ ($\mu=1,2$) are the covariant derivatives 
acting on $z:\Omega \rightarrow \field{C}^N$, defined by the formula
\begin{equation}
D_{\mu}z =\partial_{\mu}z- (z^{\dagger}\cdot \partial_{\mu}z)z\,,
\label{covader}
\end{equation}
where $\partial_{\mu}=\partial_{\xi^{\mu}}$. The action $S$ does not depend 
on the choice of a representative of the class $[z]$. 
As usual, the symbol $\dagger$ denotes Hermitian conjugation, whereas the 
Hermitian inner product of $z=(z_1, \ldots, z_N)$ and 
$w=(w_1, \ldots, w_N)$ in $\field{C}^N$ is denoted by
\begin{equation}
<z,w>=z^{\dagger}\cdot w=\sum_{j=1}^{N}\bar{z}_jw_j.
\end{equation}

Introducing
\begin{equation}
z=\frac{f}{|f|}, \qquad |f|=(f^{\dagger}\cdot f)^{\frac{1}{2}}\,,
\end{equation}
the action functional (\ref{action}) can be expressed as 
\begin{equation}
S=\frac{1}{4}\int_{\Omega}\frac{1}{f^{\dagger}\cdot f}
(\partial f^{\dagger} P \bar{\partial} f+\bar{\partial}
f^{\dagger}P \partial f) d\xi d\bar{\xi}\,, \label{action2}
\end{equation}
where $\partial$ and $\bar{\partial}$ denote the partial derivatives 
with respect to $\xi$ and $\bar{\xi}$, respectively, {i.e.,}
\begin{equation}
\partial=\frac{1}{2}\left(\partial_{\xi^1}-i\partial_{\xi^2}\right)\,, 
\qquad
\bar{\partial}=\frac{1}{2}\left(\partial_{\xi^1}+i\partial_{\xi^2}\right)\,.
\label{partials}
\end{equation}
The $N \times N$ matrix $P$ is an orthogonal projector on the orthogonal 
complement of the complex line in $\field{C}^{N}$. 
Therefore, 
\begin{equation}
P=I_{N}-\frac{1}{f^{\dagger}\cdot f} f \otimes f^{\dagger}\,,
\label{projector}
\end{equation}
where $I_{N}$ is the $N \times N$ identity matrix.
Since $P$ is an orthogonal projector it has the properties
\begin{equation}
P^{\dagger}=P\,, \qquad P^2=P\,.
\end{equation}

The map $[z]$ is determined by a solution of the Euler-Lagrange equations 
which are associated with the action (\ref{action2}).
In the homogeneous coordinates $f$, the equations of motion take the form 
of a conservation law
\begin{equation}
\partial K - \bar{\partial} K^{\dagger}=0\,, \label{conservation}
\end{equation}
where $K$ and $K^{\dagger}$ are $N \times N$ matrices given by
\begin{eqnarray}
&K=[\bar{\partial}P, P]=\frac{1}{f^{\dagger}\cdot f}
\left(\bar{\partial}f\otimes f^{\dagger}-f\otimes 
\bar{\partial}f^{\dagger}\right)
+\frac{f \otimes f^{\dagger}}{(f^{\dagger}\cdot f)^2}
\big(\bar{\partial}f^{\dagger}\cdot f - f^{\dagger}\cdot\bar{\partial}f
\big)\,, \nonumber \\
\label{defmatK} \\
&K^{\dagger}=-[{\partial}P, P]=\frac{1}{f^{\dagger}\cdot f}
\left(f\otimes \partial f^{\dagger}-\partial f\otimes 
f^{\dagger}\right)
+\frac{f \otimes f^{\dagger}}{(f^{\dagger}\cdot f)^2}
\big({\partial}f^{\dagger}\cdot f - f^{\dagger}\cdot {\partial}f
\big)\,.\nonumber 
\end{eqnarray}

Using the projector, the Euler-Lagrange equations (\ref{conservation}) 
can also be written in the form of a conservation law
\begin{equation}
\partial [\bar{\partial} P, P] + \bar{\partial} [\partial P, P]=0\,.
\label{consevation2}
\end{equation}
Through explicit calculation one can verify that the complex-valued functions
\begin{equation}
J=\frac{1}{f^{\dagger}\cdot f}\partial f^{\dagger} P \partial f\,,
\qquad
\bar{J}=\frac{1}{f^{\dagger}\cdot f}\bar{\partial} f^{\dagger} 
P \bar{\partial} f\,, \label{defJ}
\end{equation}
satisfy
\begin{equation}
\bar{\partial} J=0\,, \qquad \partial \bar{J}=0\,,
\end{equation}
for any solution $f$ of the equations of motion (\ref{conservation}).

Note that the action (\ref{action}), as well as $J$ and $\bar{J}$, 
are invariant under a global ${U}(N)$ transformation, {i.e.,}
$f \rightarrow u f$, where $u \in {U}(N)$. 
Due to this invariance, without loss of generality, we can set one 
of the components of the vector field $f$ equal to $1$. 
For instance, $f_1 = 1$. 
Consequently, the $\field{C}P^{N-1}$ model can be expressed in one less
variable through the relation
\begin{equation}
w_{i-1}=\frac{f_i}{f_1}\,, \qquad i=2, \ldots, N-1\,. \label{coordinw}
\end{equation}

\section{The Weierstrass formula for immersion 
\label{secweierstrass}} 
For a given projector $P$ satisfying the conservation law 
(\ref{consevation2}), we give the analytical description of a $2D$ 
smooth orientable surface $\mathcal{F}$ immersed in the ${su}(N)$ algebra. 
This is accomplished by constructing an exact ${su}(N)$ matrix-valued 
$1$-form $dX$ for which its ``potential,'' which is a matrix-valued 
$0$-form $X$, determines a surface immersed in the ${su}(N)$ algebra. 
Once the $0$-form $X$ is calculated, we can treat the components of $X$ 
as the coordinates of a surface in ${su}(N)$ and, hence, we can compute
an explicit formula for immersion. 
In what follows, we shall refer to this as the generalized Weierstrass 
formula for immersion. 
Next, we investigate some geometrical properties of the surface $\mathcal{F}$ 
in the ${su}(N)$ algebra.

In order to construct and investigate surfaces in multi-dimensional 
spaces by analytical methods it is convenient to identify the ${su}(N)$
algebra with the ($N^2-1$)-dimensional Euclidean space through the relation
\begin{equation}
\field{R}^{N^2-1} \simeq {su}(N)\,.  \label{isomorphism}
\end{equation}
For the sake of uniformity, we use the following definition of scalar 
product on ${su}(N)$
\begin{equation}
<A, B>=-\frac{1}{2}\rm{tr}(AB)\,,
\end{equation}
where $A,B \in su(N)$.

Let us assume that the matrix $K$ in (\ref{defmatK}) is constructed
from a solution $P$ of the Euler-Lagrange equation (\ref{consevation2})
defined on some connected and simply connected domain 
$\Omega \subset \field{C}$. 
According to Poincar\'e's lemma, there then exists a closed matrix-valued 
$1$-form, 
\begin{equation}
dX=i(K^{\dagger}d\xi + K d\bar{\xi})\,,  \label{complex1formimmer}
\end{equation}
which is also exact and takes its values in the $su(N)$ algebra of 
skew-Hermitian matrices. 
This means that $X$ is a well-defined $su(N)$ real-valued function 
on $\Omega$ and 
\begin{equation}
\partial X=i K^{\dagger}\,, \qquad \bar{\partial} X=iK\,.
\label{deftangent}
\end{equation}
It follows from the closedness of the $1$-form $dX$ that the integral
\begin{equation}
i\int_{\gamma} (K^{\dagger}d\xi + K d\bar{\xi})=X(\xi, \bar{\xi})\,,
\label{intforimm}
\end{equation}
is locally independent of the path of integration.
As a matter of fact, the integral only depends on the end points of the 
curve $\gamma$ in $\field{C}$.

The integral (\ref{intforimm}) defines a mapping
\begin{equation}
X: \Omega \ni (\xi, \bar{\xi}) \rightarrow X (\xi, \bar{\xi})
\in su(N)\,, \label{mappforimm}
\end{equation}
which is called the generalized Weierstrass formula for immersion 
\cite{Grundlandetal, GrundlandSnobl}.

As a consequence of (\ref{mappforimm}), we can determine a surface 
$\mathcal{F}$ in $su(N)$ from a solution $f$ of the Euler-Lagrange equation 
(\ref{conservation}) defined on the domain $\Omega \subset \field{C}$.

The complex tangent vectors to a surface $\mathcal{F}$ are given by
(\ref{deftangent}) using (\ref{defmatK}). 
For the components of the induced metric one gets
\begin{eqnarray}
&g_{\xi\xi}\equiv (\partial X, \partial X)=-J\,, \qquad
g_{\bar{\xi}\bar{\xi}}\equiv (\bar{\partial} X, \bar{\partial} X)
=-\bar{J}\,, \nonumber \\
&g_{\xi\bar{\xi}}=g_{\bar{\xi}\xi} \equiv (\partial X, \bar{\partial} X)
=q\,, \label{compofmetric}
\end{eqnarray}
where $J$ and $\bar{J}$ are holomorphic functions defined in (\ref{defJ}) 
and the quantity $q$ is a positive real-valued function given by
\begin{equation}
q=\frac{1}{f^{\dagger}\cdot f} \bar{\partial}f^{\dagger} P \partial 
f \geq 0 \,. \label{definedq}
\end{equation}
Thus, the first fundamental form of a surface $\mathcal{F}$ takes the form
\begin{equation}
I=-Jd\xi^2 +2q d\xi d\bar{\xi} - \bar{J} d{\bar{\xi}}^2\,. 
\label{firstfundform}
\end{equation}
Using the Schwartz inequality, it was shown in 
\cite{Grundlandetal, GrundlandSnobl} that this first 
fundamental form (\ref{firstfundform}) is positive definite. 

The scalar curvature is given by 
\begin{equation}
\mathcal{K}=\frac{1}{2\,\sqrt{g}}\bar{\partial}\left[\frac{q}{\sqrt{g}}
\partial \ln\left(-\frac{q^2}{{J}}\right)\right]\,,
\qquad \rm{if} \quad J\ne0
\label{gaussiancur1}
\end{equation}
and
\begin{equation}
\mathcal{K}=-q^{-1}\bar{\partial}\partial\ln q\,, \qquad \rm{if} 
\quad J=0\,,
\label{gaussiancur2}
\end{equation}
where
\begin{equation}
g=\rm{det}(g_{ij})=|J|^2-q^2\,, 
\end{equation}
and the indices $i$ and $j$ stand for $\xi$ and $\bar{\xi}$, respectively.

Let us now discuss the existence of certain classes of surfaces in the 
$su(N)$ algebra when the $\field{C}P^{N-1}$ equations are subjected to 
specific differential constraints (DCs). These constraints allow us to 
reduce the overdetermined system to a system admitting first integrals. 
Doing so considerably simplifies the process of solving the initial 
$\field{C}P^{N-1}$ equations (\ref{conservation}). Consequently, certain 
classes of non-splitting solutions can be constructed and they provide us 
with an explicit, simplified form of Weierstrass formula for immersion 
of a surface in $su(N)$.

\newtheorem{pre}{Proposition}
\begin{pre}
If the complex-valued vector function
\begin{equation}
\field{C}\ni \xi\rightarrow f(\xi) \in \field{C}^N \backslash \{0\}
\end{equation}
satisfies both the equations (\ref{conservation}) for the 
$\field{C}P^{N-1}$ model equations and the differential constraints
\begin{equation}
f^{\dagger}\cdot \partial f -\partial f^{\dagger} \cdot f=0\,,
\qquad
f^{\dagger}\cdot \bar{\partial} f -\bar{\partial} 
f^{\dagger} \cdot f=0\,,
\label{constraintforspecialcase}
\end{equation}
then the generalized Weierstrass formula for immersion of a 
surface $\mathcal{F}$ in the $su(N)$ algebra has the form
\begin{equation}
X(\xi, \bar{\xi})=i\int_{\gamma} \frac{f \otimes 
\partial f^{\dagger} -(\partial f^{\dagger} \cdot f) \widetilde{P}}
{f^{\dagger} \cdot f} 
d \xi + 
\frac{\bar{\partial} f \otimes 
f^{\dagger} -(f^{\dagger} \cdot \bar{\partial} f) \widetilde{P}}
{f^{\dagger} \cdot f}
d \bar{\xi}\,,
\label{generalizedweierstrassforspecial}
\end{equation}
where $\widetilde{P}=I_{N}-P$. The first fundamental form is given by
\begin{equation}
I=-{J_1}d\xi^2 +2\left(\frac{\bar{\partial} f ^{\dagger} \cdot 
\partial f}{f^{\dagger} \cdot f}- \frac{(\bar{\partial}f^{\dagger} \cdot f) 
(f^{\dagger} \cdot {\partial} f)}{(f^{\dagger} \cdot f)^2}\right) 
d\xi d\bar{\xi} - 
\bar{{{J_1}}} d{\bar{\xi}}^2\,,
\label{firstfundformforspecial}
\end{equation}
where ${J_1}$ and $\bar{{J_1}}$ are holomorphic functions,
\begin{equation}
{J_1}=\frac{\partial f^{\dagger} \cdot \partial f}{f^{\dagger} \cdot f} 
-\left(\frac{f^{\dagger} \cdot \partial f}{f^{\dagger} 
\cdot f}\right)^2\,,
\qquad
\bar{{J_1}}=\frac{\bar{\partial} f^{\dagger} \cdot \bar{\partial} f}
{f^{\dagger} 
\cdot f} -\left(\frac{\bar{\partial} f^{\dagger} \cdot f}{f^{\dagger} 
\cdot f}\right)^2\,, \label{newdefJ}
\end{equation}
which satisfy
\begin{equation}
\bar{\partial} {J_1}=0\,, \qquad \partial \bar{{J_1}}=0\,,
\end{equation}
whenever (\ref{conservation}) and (\ref{constraintforspecialcase}) hold.
\end{pre}
\vskip 5pt
\noindent
\textbf{Proof} If we append the two DCs in (\ref{constraintforspecialcase}) 
to the $\field{C}P^{N-1}$ equations (\ref{conservation}) then the matrices 
$K$ and $K^{\dagger}$ in (\ref{defmatK}), become
\begin{eqnarray}
{K_1}=\frac{1}{f^{\dagger} \cdot f} (\bar{\partial} f \otimes 
f^{\dagger} -f \otimes \bar{\partial} f^{\dagger})\,, \nonumber \\
K_1^{\dagger}=\frac{1}{f^{\dagger} \cdot f} (f \otimes \partial 
f^{\dagger}-\partial f \otimes f^{\dagger})\,. \label{newdefK}
\end{eqnarray}
Hence, the Weierstrass formula for immersion takes the form 
\begin{eqnarray}
X(\xi, \bar{\xi})&=&i\int_{\gamma} (K_1^{\dagger}d\xi + 
{K_1} d\bar{\xi}) \nonumber \\
&=& i\int_{\gamma}\frac{f \otimes \partial 
f^{\dagger}-\partial f \otimes f^{\dagger}}{f^{\dagger} \cdot f} d\xi
+
\frac{\bar{\partial} f \otimes 
f^{\dagger} -f \otimes \bar{\partial} f^{\dagger}}{f^{\dagger} \cdot f}
d \bar{\xi}\,.
\end{eqnarray}
On the other hand, from (\ref{conservation}), we are able to deduce 
that the matrix $K$ can be decomposed as
\begin{equation}
K=M+L\,, \label{decopositionofK}
\end{equation}
where
\begin{equation}
M=(I-P)\bar{\partial} P\,, \qquad L=-\bar{\partial} P (I-P)\,.
\label{mandlwithp}
\end{equation}
It can be shown that the matrices $M$ and $L$ satisfy the same 
conservation laws
(\ref{conservation}) as the matrix $K$, {e.g.,} 
\begin{equation}
\partial M= \bar{\partial} M^{\dagger}\,, \qquad \partial L=
\bar{\partial} L^{\dagger}\,. \label{conservationlawsforML}
\end{equation}
Note that the two conservation laws in (\ref{conservationlawsforML}) 
are not independent since $M$ and $L$ differ by a total divergence, 
\begin{equation}
M=L+\bar{\partial} P\,. \label{totaldivergence}
\end{equation}
Taking into account the overdetermined system composed of the conservation
laws (\ref{conservation}) and DCs (\ref{constraintforspecialcase}) 
for the function $f$, the matrices $M$ and $L$ become
\begin{eqnarray}
&&{M_1}=-\frac{f \otimes 
\bar{\partial} f^{\dagger} -(f^{\dagger}\cdot \bar{\partial} f) \widetilde{P}} 
{f^{\dagger} \cdot f} \,,
\qquad
M_1^{\dagger}=-\frac{{\partial} f 
\otimes 
f^{\dagger} -(\partial f^{\dagger} \cdot f) \widetilde{P}}
{f^{\dagger} \cdot f} \,,
\nonumber \\
&&{L_1}=\frac{\bar{\partial} f \otimes 
f^{\dagger} -(f^{\dagger} \cdot \bar{\partial} f) \widetilde{P}}
{f^{\dagger} \cdot f} \,,
\qquad
L_1^{\dagger}=\frac{f \otimes 
\partial f^{\dagger} -(\partial f^{\dagger} \cdot f) \widetilde{P}}
{f^{\dagger} \cdot f} \,.
\label{newMandL}
\end{eqnarray}
As a consequence of the conservation laws (\ref{conservationlawsforML}) 
for the matrices ${M_1}$ and ${L_1}$, the Weierstrass formula for 
immersion (\ref{intforimm}) takes the following simple form
\begin{eqnarray}
X(\xi, \bar{\xi})&=&i\int_{\gamma} (M_1^{\dagger}d\xi + 
{M_1} d\bar{\xi}) \nonumber \\
&=&-i\int_{\gamma} \frac{{\partial} f 
\otimes 
f^{\dagger} -(\partial f^{\dagger} \cdot f) \widetilde{P}}
{f^{\dagger} \cdot f} d\xi
+ \frac{f \otimes 
\bar{\partial} f^{\dagger} -(f^{\dagger}\cdot \bar{\partial} f) \widetilde{P}}
{f^{\dagger} \cdot f}
d \bar{\xi}\,, \label{immersionforM}
\end{eqnarray}
or
\begin{eqnarray}
X(\xi, \bar{\xi})&=&i\int_{\gamma} (L_1^{\dagger}d\xi + 
{L_1} d\bar{\xi}) \nonumber \\
&=& i\int_{\gamma} \frac{f \otimes 
\partial f^{\dagger} -(\partial f^{\dagger} \cdot f) \widetilde{P}}
{f^{\dagger} \cdot f} d\xi +
\frac{\bar{\partial} f \otimes 
f^{\dagger} -(f^{\dagger} \cdot \bar{\partial} f) \widetilde{P}}
{f^{\dagger} \cdot f}
d \bar{\xi}\,, \label{immersionforL}
\end{eqnarray}
respectively. 
As a consequence of 
(\ref{totaldivergence}), (\ref{immersionforM}) and (\ref{immersionforL}),  
it can be shown that the two different Weierstrass data $({L_1}, 
L_1^{\dagger})$ or $({M_1}, 
M_1^{\dagger})$ correspond to different parametrizations 
of 
the same surface $\mathcal{F}$ in the $su(N)$ algebra. 

In this case, the quantity $J$ takes the simple form 
\begin{equation}
{J_1}=\frac{\partial f^{\dagger} \cdot \partial f}{f^{\dagger} \cdot f} 
-\left(\frac{f^{\dagger} \cdot \partial f}{f^{\dagger} 
\cdot f}\right)^2\,,
\qquad
\bar{{J_1}}=\frac{\bar{\partial} f^{\dagger} \cdot \bar{\partial} f}
{f^{\dagger} 
\cdot f} -\left(\frac{\bar{\partial} f^{\dagger} \cdot f}{f^{\dagger} 
\cdot f}\right)^2\,. \label{newdefJ2}
\end{equation}
Using the conservation laws (\ref{conservation}) and DCs 
(\ref{constraintforspecialcase}) for the function $f$, we find that
${J_1}$ is a holomorphic function, {e.g.,} $\bar{\partial} J_1=0$ 
whenever (\ref{conservation}) and (\ref{constraintforspecialcase}) hold. 
As a consequence of (\ref{immersionforM}), 
(\ref{immersionforL}) and (\ref{newdefJ2}), the components of the 
induced metric are
\begin{equation}
g_{\xi\xi}=-{J_1}\,, \qquad g_{\bar{\xi}\bar{\xi}}=-\bar{{J_1}}\,, 
\qquad g_{\xi\bar{\xi}}=\frac{\bar{\partial} f ^{\dagger} \cdot 
\partial f}{f^{\dagger} \cdot f}- \frac{(\bar{\partial}f^{\dagger} \cdot f) 
(f^{\dagger} \cdot {\partial} f)}{(f^{\dagger} \cdot f)^2}\,, 
\label{inducedmetspecial}
\end{equation}
which completes the proof. $\;\;\;\;\hspace{7cm}\square$ 

Note that the complex-valued vector function $\field{C} \ni \xi 
\rightarrow f(\xi) \in \field{C}^N \backslash \{0\}$ is a 
holomorphic $(\bar{\partial} f=0)$ solution of the $\field{C}P^{N-1}$ 
model (\ref{conservation}) if and only if the generalized Weierstrass 
formula for the immersion of a surface $\mathcal{F}$ has the 
skew-Hermitian form
\begin{equation}
X(\xi, \bar{\xi})=-i P \in su(N)\,. \label{skewHermitian}
\end{equation}
If $f$ is holomorphic, {\it i.e.,} $\bar{\partial}f=0$, then by virtue of 
equations (\ref{mandlwithp}) and the differential consequences of the 
identity $(I_{N}-P)P=0$, we obtain
\begin{equation}
M=0\,, \qquad \bar{\partial}P\,P=0\,. \label{newdefmwithp}
\end{equation}
Using the differential consequences for the projector $P$, we get 
\begin{eqnarray}
\bar{\partial}P\,P=0\,, \qquad P\,\partial P=0\,, \nonumber \\
\bar{\partial}P=P\,\bar{\partial}P\,, \qquad \partial P= \partial P P\,.
\label{diffconsforp}
\end{eqnarray}
Substituting (\ref{diffconsforp}) into (\ref{defmatK}), we obtain
\begin{equation}
K=-\bar{\partial} P\,, \qquad K^{\dagger}=-\partial P\,.
\label{newkintermsofp}
\end{equation}
Hence, the Weierstrass formula for immersion (\ref{intforimm}) of 
$\mathcal{F}$ is expressed in terms of the projector $P$ and is a 
skew-Hermitian matrix given by (\ref{skewHermitian}). 
This result coincides with the one obtained in \cite{Hussin}.

The converse is also true. 
Indeed, if we assume that the Weierstrass formula for immersion 
(\ref{intforimm}) of $\mathcal{F}$ is a projector $P$ then the 
differential of $X$ leads to (\ref{newkintermsofp}). 
Using the differential consequences of the relation $P^2=P$, we obtain the 
relations (\ref{diffconsforp}) which lead to $M=0$. 
In view of equations (\ref{mandlwithp}), this implies that, 
in the generic case, solutions of the $\field{C}P^{N-1}$ model 
(\ref{conservation}) must be holomorphic.

Also, note that in the case of the holomorphic solutions of the 
$\field{C}P^{N-1}$ model, the corresponding complex-valued function 
(\ref{defJ}) vanishes, {i.e.,}
\begin{equation}
J=\frac{1}{f^{\dagger}\cdot f}\partial f^{\dagger} P \partial f =0\,.
\end{equation}

An analogous statement can be made for anti-holomorphic solutions 
(${\partial} f=0$) of equation (\ref{conservation}). 
For this case, we have
\begin{equation}
L=0\,, \qquad P\bar{\partial}P=0\,, \qquad \partial P\,P=0\,.
\end{equation}
Hence, from (\ref{defmatK}), the matrices $K$ and $K^{\dagger}$ become 
\begin{equation}
K=\bar{\partial} P\,, \qquad K^{\dagger}=\partial P\,.
\end{equation}
Finally, one can see that the Weierstrass formula for immersion of 
$\mathcal{F}$ is the skew-Hermitian form
\begin{equation}
X(\xi, \bar{\xi})=i P \in su(N)\,.
\end{equation}

\section{The Lie-point symmetries of the $\field{C}P^{N-1}$ sigma models 
\label{secliepointsymm}}
In this section, we present the explicit formulae for the Lie-point 
symmetries of the $\field{C}P^{N-1}$ model (\ref{conservation}) for
arbitrary $N$. 
To do so, we first compute the symmetries for the $\field{C}P^{1}$, 
$\field{C}P^{2}$ and $\field{C}P^{3}$ models.
We then generalize the results to the $\field{C}P^{N-1}$ case by induction. 
For the computation of the Lie-point symmetries, we search for the most 
general (point) transformations of the independent and dependent variables 
which leave the solution set of (\ref{conservation}) invariant. 
Locally, such transformations are given by a vector field of the form 
\cite{Olver}
\begin{equation}
\vec{v}=\eta^1\partial+\eta^2\bar{\partial}+\sum_{j=1}^{N-1}\Phi_j^1
\partial_{w_j}+\sum_{j=1}^{N-1}\Phi_j^2\partial_{\bar{w}_j}\,,
\label{genvector}
\end{equation}
where $\eta^1$, $\eta^2$, $\Phi_j^1$ and $\Phi_j^2$ are functions of
$\xi$, $\bar{\xi}$ and the affine coordinates $w_1, \bar{w}_1,$ $\ldots$,
$w_{N-1}, \bar{w}_{N-1}$. 
According to the symmetry criterion \cite{Olver}, the second prolongation 
of $\vec{v}$ acting on (\ref{conservation}) must vanish on the solution 
set of (\ref{conservation}). 
This requirement leads to the so-called determining equations, 
whose solution yields the functions 
$\eta^1$, $\eta^2$, $\Phi_j^1$ and $\Phi_j^2$. 

Generating the determining equations is entirely algorithmic.
Reducing and solving them can be done fully automatic with sophisticated 
software, or, perhaps more reliably, by interactively adding information 
extracted from the simplest determining equations before computing 
the full set. 
Many software packages have been written to perform Lie symmetry computations.
In-depth reviews of such packages can be found in 
\cite{Champagneetal, HeremanCRC, HeremanMCM, ButcherCPC}.

For low dimensions, {e.g.,} for $N\le4,$ we did the computations 
independently with {\sc SYMMGRP.MAX} and by hand. 
For the latter, we took advantage of the discrete symmetries of the model.
For the $\field{C}P^{N-1}$ models with $N\geq4$, after eliminating all 
single-term determining equations and their differential consequences, 
we were left with several hundred of determining equations.
Using {\sc SYMMGRP.MAX} interactively, these determining equations were 
further reduced and eventually completely solved.

We now discuss the Lie-point symmetries of the $\field{C}P^{N-1}$ models 
for $N=2, 3$, and $4$, separately. 

The equations for the $\field{C}P^1$ model, expressed in terms of 
the homogeneous coordinate $w_1$ defined in (\ref{coordinw}), are given by
\begin{eqnarray}
\partial\bar{\partial}w_1-\frac{2\bar{w}_1}{A_1}\partial w_1 
\bar{\partial}w_1=0\,,   \qquad \qquad 
\partial\bar{\partial}\bar{w}_1-\frac{2{w_1}}{A_1}\partial 
\bar{w}_1 \bar{\partial}\bar{w}_1=0\,, \label{cp1}
\end{eqnarray} 
where $A_1=1+w_1\bar{w}_1$. The general solution of the determining equations 
associated with 
vector field (\ref{genvector}) is given by
\begin{eqnarray}
&&\eta^1=\eta^1(\xi)\,,\qquad \eta^2=\eta^2(\bar{\xi})\,,\nonumber \\
&&\Phi_1^1=\alpha_1 {w_1}^{\!\!2}+\beta_1 w_1+\gamma_1\,,\nonumber\\
&&\Phi_1^2=\gamma_1 {\bar{w}_1}^{\,2}-\beta_1 \bar{w}_1+\alpha_1\,,
\end{eqnarray}
where $\eta^1$ and $\eta^2$ are arbitrary functions of $\xi$ and 
$\bar{\xi}$, respectively and $\alpha_1$, $\beta_1$ and $\gamma_1$ are 
arbitrary constants. Thus, the corresponding symmetry algebra 
$\mathcal{L}_1$ is spanned by five generators, namely 
\begin{eqnarray}
&&X_1=\eta^1(\xi)\partial\,,\qquad X_2=\eta^2(\bar{\xi})\bar{\partial}
\,, \nonumber \\
&&X_3={w_1}{\!\!^2}\partial_{w_1}+\partial_{\bar{w}_1}\,,\nonumber \\
&&X_4=w_1 \partial_{w_1}-\bar{w}_1\partial_{\bar{w}_1}\,,\nonumber \\
&&X_5=\partial_{w_1}+{\bar{w}_1}^{\,2}\partial_{\bar{w}_1}\,. 
\label{cp1symm}
\end{eqnarray} 
The algebra $\mathcal{L}_1$ can be decomposed as a direct sum of two
infinite-dimensional simple Lie algebras and the $su(2)$ algebra 
generated by $\{X_3,X_4,X_5\}$, {i.e.,} 
\begin{equation}
\mathcal{L}_1=\{X_1\} \oplus \{X_2\} \oplus su(2)\,.
\end{equation}

Likewise, in terms of homogeneous coordinates $w_1$ and $w_2$ in 
(\ref{coordinw}), the equations for the $\field{C}P^2$ model read
\begin{eqnarray}
&&\partial\bar{\partial}w_1-\frac{2\bar{w}_1}{A_2}\partial w_1 
\bar{\partial} w_1-\frac{\bar{w}_2}{A_2}(\partial w_1 \bar{\partial} 
w_2 + \bar{\partial} w_1 \partial w_2)=0\,,\nonumber\\
&&\partial\bar{\partial}w_2-\frac{2\bar{w}_2}{A_2}\partial w_2 
\bar{\partial} w_2-\frac{\bar{w}_1}{A_2}(\partial w_1 \bar{\partial} 
w_2 + \bar{\partial} w_1 \partial w_2)=0\,,\nonumber\\
&&\partial\bar{\partial}\bar{w}_1-\frac{2{w}_1}{A_2}\partial 
\bar{w}_1 \bar{\partial} \bar{w}_1-\frac{{w}_2}{A_2}(\bar{\partial} 
\bar{w}_1 {\partial} \bar{w}_2 + {\partial} \bar{w}_1 \bar{\partial} 
\bar{w}_2)=0\,,\nonumber\\
&&\partial\bar{\partial}\bar{w}_2-\frac{2{w}_2}{A_2}\partial \bar{w}_2 
\bar{\partial} \bar{w}_2-\frac{{w}_1}{A_2}(\bar{\partial} 
\bar{w}_1 {\partial} \bar{w}_2 + {\partial} \bar{w}_1 
\bar{\partial} \bar{w}_2)=0\,, \label{cp2}
\end{eqnarray}
where $A_2=1+w_1\bar{w}_1+w_2\bar{w}_2$. 
Upon integration, the determining equations yield 
\begin{eqnarray}
&&\eta^1=\eta^1(\xi)\,,\qquad \eta^2=\eta^2(\bar{\xi})\,,\nonumber \\
&&\Phi_1^1=k_1{w_1}^{\!\!2}+k_2 w_1 w_2+k_4 w_1 +k_5 w_2 +k_6\,,
\nonumber\\
&&\Phi_2^1=k_2{w_2}^{\!\!2}+k_1 w_1 w_2+k_3 w_2 +k_7 w_1 +k_8\,,
\nonumber\\
&&\Phi_1^2=k_6{\bar{w}_1}^{\,2}+k_8 \bar{w}_1 \bar{w}_2-k_4 \bar{w}_1 -k_7 
\bar{w}_2 +k_1\,,\nonumber\\
&&\Phi_2^2=k_8{\bar{w}_2}^{\,2}+k_6 \bar{w}_1 \bar{w}_2-k_3 \bar{w}_2 -k_5 
\bar{w}_1 +k_2\,,
\end{eqnarray}
where $k_i$ ($i=1,\ldots,8$) are arbitrary constants. 
The associated symmetry algebra $\mathcal{L}_2$ of (\ref{cp2}) is 
thus spanned by the following ten generators:
\begin{eqnarray}
&&X_1=\eta^1(\xi)\partial\,,\qquad X_2=\eta^2(\bar{\xi})\bar{\partial}
\,, \nonumber \\
&&X_3={w_1}^{\!\!2}\partial_{w_1}+w_1w_2\partial_{w_2}+\partial_{\bar{w}_1}
\,, \nonumber \\
&&X_4=w_1w_2\partial_{w_1}+{w_2}^{\!\!2}\partial_{w_2}+\partial_{\bar{w}_2}
\,, \nonumber \\
&&X_5=w_2\partial_{w_2}-\bar{w}_2\partial_{\bar{w}_2}\,, \nonumber \\
&&X_6=w_1\partial_{w_1}-\bar{w}_1\partial_{\bar{w}_1}\,, \nonumber \\
&&X_7=w_2\partial_{w_1}-\bar{w}_1\partial_{\bar{w}_2}\,, \nonumber \\
&&X_8=\partial_{w_1}+{\bar{w}_1}^{\,2}\partial_{\bar{w}_1}+\bar{w}_1
\bar{w}_2\partial_{\bar{w}_2}\,,\nonumber\\
&&X_9=w_1\partial_{w_2}-\bar{w}_2\partial_{\bar{w}_1}\,, \nonumber \\
&&X_{10}=\partial_{w_2}+\bar{w}_1\bar{w}_2\partial_{\bar{w}_1}+
{\bar{w}_2}^{\,2}\partial_{\bar{w}_2}\,. \label{cp2symm}
\end{eqnarray}
As in the previous case, the symmetry algebra $\mathcal{L}_2$ can be 
decomposed as a direct sum of two infinite-dimensional simple Lie algebras 
and the $su(3)$ algebra.

In like fashion, in terms of $w_1$, $w_2$ and $w_3$ in (\ref{coordinw}), 
the equations for the $\field{C}P^3$ model are
\begin{eqnarray}
&&\!\!\!\!\!\!\!\!\!\!
\partial\bar{\partial}w_1-\frac{2\bar{w}_1}{A_3}\partial w_1 
\bar{\partial} w_1-\frac{\bar{w}_2}{A_3}(\partial w_1 \bar{\partial} w_2 
+ \bar{\partial} w_1 \partial w_2)-\frac{\bar{w}_3}{A_3}(\partial w_1 
\bar{\partial} w_3 + \bar{\partial} w_1 \partial w_3)=0\,,\nonumber\\
&&\!\!\!\!\!\!\!\!\!\!
\partial\bar{\partial}w_2-\frac{2\bar{w}_2}{A_3}\partial w_2 
\bar{\partial} w_2-\frac{\bar{w}_1}{A_3}(\partial w_1 \bar{\partial} w_2 
+ \bar{\partial} w_1 \partial w_2)-\frac{\bar{w}_3}{A_3}(\partial w_2 
\bar{\partial} w_3 + \bar{\partial} w_2 \partial w_3)=0\,,\nonumber\\
&&\!\!\!\!\!\!\!\!\!\!
\partial\bar{\partial}w_3-\frac{2\bar{w}_3}{A_3}\partial w_3 
\bar{\partial} w_3-\frac{\bar{w}_1}{A_3}(\partial w_1 \bar{\partial} w_3 
+ \bar{\partial} w_1 \partial w_3)-\frac{\bar{w}_2}{A_3}(\partial w_2 
\bar{\partial} w_3 + \bar{\partial} w_2 \partial w_3)=0\,,\nonumber\\
&&\!\!\!\!\!\!\!\!\!\!
\partial\bar{\partial}\bar{w}_1-\frac{2 w_1}{A_3}\partial \bar{w}_1 
\bar{\partial} \bar{w}_1-\frac{w_2}{A_3}(\partial \bar{w}_1 
\bar{\partial} \bar{w}_2 + \bar{\partial} \bar{w}_1 \partial \bar{w}_2)-
\frac{w_3}{A_3}(\partial \bar{w}_1 \bar{\partial} \bar{w}_3 + 
\bar{\partial} \bar{w}_1 \partial \bar{w}_3)=0\,,\nonumber\\
&&\!\!\!\!\!\!\!\!\!\!
\partial\bar{\partial}\bar{w}_2-\frac{2 w_2}{A_3}\partial 
\bar{w}_2 \bar{\partial} \bar{w}_2-\frac{w_1}{A_3}(\partial 
\bar{w}_1 \bar{\partial} \bar{w}_2 + \bar{\partial} \bar{w}_1 \partial 
\bar{w}_2)-\frac{w_3}{A_3}(\partial \bar{w}_2 \bar{\partial} \bar{w}_3 
+ \bar{\partial} \bar{w}_2 \partial \bar{w}_3)=0\,,\nonumber\\
&&\!\!\!\!\!\!\!\!\!\!
\partial\bar{\partial}\bar{w}_3-\frac{2 w_3}{A_3}\partial 
\bar{w}_3 \bar{\partial} \bar{w}_3-\frac{w_1}{A_3}(\partial 
\bar{w}_1 \bar{\partial} \bar{w}_3 + \bar{\partial} \bar{w}_1 
\partial \bar{w}_3)-\frac{w_2}{A_3}(\partial \bar{w}_2 \bar{\partial} 
\bar{w}_3 + \bar{\partial} \bar{w}_2 \partial \bar{w}_3)=0\,,
\nonumber\\ \label{cp3}
\end{eqnarray}
where $A_3=1+w_1\bar{w}_1+w_2\bar{w}_2+w_3\bar{w}_3$. 
After straightforward but long calculations the determining equations yield 
\begin{eqnarray}
&&\eta^1=\eta^1(\xi)\,,\qquad \eta^2=\eta^2(\bar{\xi})\,,\nonumber \\
&&\Phi_1^1=c_1{w_1}^{\!\!2}+c_2 w_1 w_2+c_3 w_1 w_3 +c_7 w_1+c_{10} w_2 +
c_{11} w_3 +c_4\,,\nonumber\\
&&\Phi_2^1=c_2{w_2}^{\!\!2}+c_1 w_1 w_2+c_3 w_2 w_3 +c_{13} w_1+c_{8} w_2 
+c_{12} w_3 +c_5\,,\nonumber\\
&&\Phi_3^1=c_3{w_3}^{\!\!2}+c_1 w_1 w_3+c_2 w_2 w_3 +c_{14} w_1+c_{15} w_2 
+c_{9} w_3 +c_6\,,\nonumber\\
&&\Phi_1^2=c_4{\bar{w}_1}^{\,2}+c_5 \bar{w}_1 \bar{w}_2+c_6 \bar{w}_1 
\bar{w}_3-c_7 \bar{w}_1 -c_{13} \bar{w}_2-c_{14} \bar{w}_3+c_1\,,
\nonumber\\
&&\Phi_2^2=c_5{\bar{w}_2}^{\,2}+c_4 \bar{w}_1 \bar{w}_2+c_6 \bar{w}_2 
\bar{w}_3-c_{10} \bar{w}_1 -c_{8} \bar{w}_2-c_{15} \bar{w}_3+c_2\,,
\nonumber\\
&&\Phi_3^2=c_6{\bar{w}_3}^{\,2}+c_4 \bar{w}_1 \bar{w}_3+c_5 \bar{w}_2 
\bar{w}_3-c_{11} \bar{w}_1 -c_{12} \bar{w}_2-c_{9} \bar{w}_3+c_3\,,
\end{eqnarray}
where $c_i$ ($i=1,\ldots,15$) are arbitrary constants. 
Hence, the generators corresponding to the symmetry algebra 
$\mathcal{L}_3$ of (\ref{cp3}) are given by
\begin{eqnarray}
&&X_1=\eta^1(\xi)\partial\,,\qquad X_2=\eta^2(\bar{\xi})\bar{\partial}
\,, \nonumber \\
&&S_i=w_i \partial_{w_i}-\bar{w}_i\partial_{\bar{w}_i}\,,
\nonumber \\
&&T_{ij}=w_i \partial_{w_j}-\bar{w}_j\partial_{\bar{w}_i}\,,
\qquad i\ne j\,,\nonumber \\
&&Y_i=w_i^2\partial_{w_i} + \sum_{j\ne i}^{3} w_iw_j\partial_{w_j}+
\partial_{\bar{w}_i}\,,\nonumber \\
&&Z_i=\bar{w}_i^2\partial_{\bar{w}_i} + \sum_{j\ne i}^{3} 
\bar{w}_i\bar{w}_j\partial_{\bar{w}_j}
+\partial_{w_i}\,,
\label{cp3symm}
\end{eqnarray}
where $i,j=1,2,3$. 
From $S_i$, $Y_i$ and $Z_i$ we get nine generators; 
from $T_{ij}$ we obtain six generators. 
The symmetry algebra $\mathcal{L}_3$ can be written as a direct sum 
of two infinite-dimensional simple Lie algebras and $su(4)$. 
The results for the low-dimensional cases reveal an emerging pattern:
the symmetry algebra is a direct sum of two infinite-dimensional Lie 
algebras and a finite-dimensional one.
Furthermore, the finite-dimensional part of the symmetry algebras for 
the $\field{C}P^1$, $\field{C}P^2$ and $\field{C}P^3$ models are 
associated with the $su(2)$, $su(3)$ and $su(4)$ algebras, respectively. 

We now turn to the $\field{C}P^{N-1}$ model for arbitrary $N$.  
In homogeneous coordinates $w_i,$ the equations are 
\begin{eqnarray}
\partial\bar{\partial}w_i-\frac{2\bar{w}_i}{A_{N-1}}\partial w_i 
\bar{\partial} w_i-
\frac{1}{A_{N-1}}\sum_{j\ne i}^{N-1} \bar{w}_j (\partial w_i 
\bar{\partial} w_j + \bar{\partial} w_i \partial w_j)=0\,,\nonumber\\
\partial\bar{\partial}\bar{w}_i-\frac{2 w_i}{A_{N-1}}\partial 
\bar{w}_i \bar{\partial} \bar{w}_i-
\frac{1}{A_{N-1}} \sum_{j\ne i}^{N-1} w_j (\partial 
\bar{w}_i \bar{\partial} \bar{w}_j + \bar{\partial} \bar{w}_i 
\partial \bar{w}_j)=0\,, \label{cpn}
\end{eqnarray}
where $i=1,2,...,N-1$ and $A_{N-1}=1+\sum_i^{N-1}w_i\bar{w}_i$. 

By induction, it can be shown that the symmetry algebra
$\mathcal{L}_{N-1}$ of (\ref{cpn}) is generated by
\begin{eqnarray}
&&X_1=\eta^1(\xi)\partial\,,\qquad X_2=\eta^2(\bar{\xi})\bar{\partial}
\,, \nonumber \\
&&S_i=w_i \partial_{w_i}-\bar{w}_i\partial_{\bar{w}_i}\,, 
\nonumber \\
&&T_{ij}=w_i \partial_{w_j}-\bar{w}_j\partial_{\bar{w}_i}\,, 
\qquad i\ne j\,,\nonumber \\
&&Y_i=w_i^2\partial_{w_i} +
\sum_{j\ne i}^{N-1} w_iw_j\partial_{w_j}+\partial_{\bar{w}_i}\,,
\nonumber \\
&&Z_i=\bar{w}_i^2\partial_{\bar{w}_i} +
\sum_{j\ne i}^{N-1}
\bar{w}_i
\bar{w}_j\partial_{\bar{w}_j}+\partial_{w_i}\,,
\label{cpNsymm}
\end{eqnarray}
where $i,j=1,2,...,N-1$. 
Furthermore, it can be shown that the symmetry algebra $\mathcal{L}_{N-1}$ 
is a direct sum of two infinite-dimensional Lie algebras and the 
$su(N)$ algebra, {i.e.,}
\begin{equation}
\mathcal{L}_{N-1}=\{X_1\} \oplus \{X_2\} \oplus su(N)\,.
\end{equation}

Finally, we consider two limiting cases:
\begin{enumerate}
 \item If $w_{N-1} \rightarrow 0$ then the $\field{C}P^{N-1}$ model 
reduces to the $\field{C}P^{N-2}$ model. 
Also, if all $N-2$ homogeneous coordinates vanish, then the 
$\field{C}P^{N-1}$ model reduces to the $\field{C}P^1$ model.
 \item If $w_i \rightarrow \frac{w}{\sqrt{N-1}}$ for $i=1,\ldots,N-1$, 
then the $\field{C}P^{N-1}$ model reduces to the $\field{C}P^1$ model.
\end{enumerate}
Hence, in the $\field{C}P^1$ case, we have a significant simplification.

\section{Immersion of surfaces into the $su(3)$ algebra arising from the 
$\field{C}P^2$ sigma model \label{secimmsu3cp2}}
In this section we explore the metric aspects of surfaces immersed in 
the $su(3)$ algebra associated with the $\field{C}P^2$ model. 
From the properties of the Hermitian matrix $\partial K$ we 
determine explicitly a moving frame on a conformally
parametrized surface $\mathcal{F}$ in $\field{R}^8$. 
We also derive the corresponding Gauss-Weingarten equations expressed 
in terms of any holomorphic solution of the $\field{C}P^2$ model. 
This investigation is a follow-up to earlier work
\cite{Grundlandetal, GrundlandSnobl}.
It allows us to communicate our new insights into the subject, as well as 
to present additional geometric characteristics of surfaces obtained 
from the model.

The assumption that the set $\{w_1, w_2\}$ is a holomorphic solution
of the equations for the $\field{C}P^2$ model 
implies that the quantity $J$ in (\ref{defJ}) vanishes.  
The induced metric on $\mathcal{F}$ given in (\ref{firstfundform}) 
is then conformal. 
In the $\field{C}P^2$ case, the $3\times 3$ projector matrix in 
(\ref{projector}) reads
\begin{equation}
P=I_3-\frac{1}{A_2} \left(\begin{array}{ccc}
1 & w_1 & w_2 \\
\bar{w}_1 & w_1\bar{w}_1 & w_2\bar{w}_1 \\
\bar{w}_2 & w_1\bar{w}_2 & w_2\bar{w}_2
\end{array}
\right)\,, \label{projectorcp2}
\end{equation}
where $I_3$ is the $3\times 3$ identity matrix. 
Assume that we are dealing with the generic case. 
That is, where the projector $P$ is a solution of the Euler-Lagrange equations 
(\ref{cp2}) such that the induced metric $g$ has a non-vanishing determinant
in some neighbourhood of a regular point 
$(\xi_0, \bar{\xi}_0) \in \Omega \subset \field{C}$. 
Further assume that a conformally parametrized surface $\mathcal{F}$, 
given by (\ref{intforimm}) and associated with the $\field{C}P^2$ model 
is described by a moving frame on $\mathcal{F}$ in $\field{R}^8$
\begin{equation}
\vec{\tau}=(\eta_1=\partial X, \eta_2=\bar{\partial} X, \eta_3, 
\ldots, \eta_8)^T\,, \label{movframe}
\end{equation}
where superscript $T$ stands for transpose. 
Here, the vectors $\eta_1, \ldots, \eta_8$ are identified with 
$3\times 3$ skew-Hermitian matrices through the isomorphism 
(\ref{isomorphism}). 
Furthermore, assume that the vectors form an orthonormal set, 
\begin{equation}
(\eta_j, \eta_k)=\delta_{jk}\,, \qquad  j,k=1,\dots, 8\,, 
\label{normaconditionswitdelta}
\end{equation}
where $\delta_{jk}$ is the Kronecker delta. 
Due to the normalization of the $su(3)$-valued function $X$ on $\Omega$, 
we can express the moving frame in (\ref{movframe}) on $\mathcal{F}$ 
in terms of
the adjoint $SU(3)$ representation. In the neighbourhood of a regular 
point $p=(\xi_0, \bar{\xi}_0) \in \field{C}$ an orthonormal moving 
frame $\vec{\tau}$ on $\mathcal{F}$ satisfies
\begin{eqnarray}
&&\eta_1=i e^{\frac{u}{2}}\phi^{\dagger} y_-\phi\,, \nonumber \\
&&\eta_2=i e^{\frac{u}{2}}\phi^{\dagger} y_{+}\phi\,, \nonumber \\
&&\eta_j=\phi^{\dagger} s_j\phi\,, \qquad j=3,\ldots, 8\,, 
\label{framecp2}
\end{eqnarray}
where $u$ is a real-valued function of $\xi$ and $\bar{\xi}$. 
The function $\phi$ in (\ref{framecp2}) belongs to $SU(3)$ and 
can be decomposed into the product of three $SU(2)$ factors, {i.e.,}
\begin{equation}
\phi=\left(\begin{array}{ccc}
1 & 0 & 0 \\
0 & a_1 & b_1 \\
0 & -\bar{b}_1 & \bar{a}_1 
\end{array}
\right)
\left(\begin{array}{ccc}
e^{i\varphi}\cos\alpha & -\sin\alpha & 0 \\
\sin\alpha & e^{-i\varphi}\cos\alpha  & 0 \\
0 & 0 & 1 
\end{array}
\right)
\left(\begin{array}{ccc}
1 & 0 & 0 \\
0 & a_2 & b_2 \\
0 & -\bar{b}_2 & \bar{a}_2 
\end{array}
\right)\,,
\label{generalsu3}
\end{equation}
where $a_i$, $b_i$ $i=1,2$ are complex-valued functions of 
$\xi$ and $\bar{\xi}$, subject to the constraints 
$|a_i|^2+|b_i|^2=1$ and $\alpha$, $\varphi$ are real-valued
functions of $\xi$, $\bar{\xi}$ $\in$ $\field{C}$. Here, the
set $\{s_1,\ldots,s_8\}$ forms an orthonormal basis of the Lie
algebra $su(3)$ ({\it {e.g.,}} the so-called Gell-Mann matrices 
\cite{Helgason})
given by
\begin{eqnarray}
&s_1=\left(\begin{array}{ccc}
0 & 0 & 0 \\
0 & 0 & -i \\
0 & -i & 0
\end{array}
\right)\,,\quad
s_2=\left(\begin{array}{ccc}
0 & 0 & 0 \\
0 & 0 & -1 \\
0 & 1 & 0
\end{array}
\right)\,,\quad
s_3=\left(\begin{array}{ccc}
0 & 0 & 0 \\
0 & -i & 0 \\
0 & 0 & i
\end{array}
\right)\,, \nonumber \\
&s_4=\frac{1}{\sqrt{3}}\left(\begin{array}{ccc}
-2i & 0 & 0 \\
0 & i & 0 \\
0 & 0 & i
\end{array}
\right)\,, 
s_5=\left(\begin{array}{ccc}
0 & -1 & 0 \\
1 & 0 & 0\\
0 & 0 & 0
\end{array}
\right)\,,
s_6=\left(\begin{array}{ccc}
0 & 0 & -1 \\
0 & 0 & 0\\
1 & 0 & 0
\end{array}
\right)\,, \nonumber \\
&s_7=\left(\begin{array}{ccc}
0 & i & 0 \\
i & 0 & 0\\
0 & 0 & 0
\end{array}
\right)\,,\qquad
s_8=\left(\begin{array}{ccc}
0 & 0 & i \\
0 & 0 & 0\\
i & 0 & 0
\end{array}
\right)\,.
\label{matricess8}
\end{eqnarray}
These matrices satisfy the following trace condition
\begin{equation}
\rm{tr}(s_i\,s_j)=-2\delta_{ij}\,.
\end{equation}
We also introduced the following notation
\begin{equation}
y_-=\frac{i}{2}(s_1-is_2)=
\left(\begin{array}{ccc}
0 & 0 & 0 \\
0 & 0 & 0\\
0 & 1 & 0
\end{array}
\right)\,, \qquad
y_+=\frac{i}{2}(s_1+is_2)=
\left(\begin{array}{ccc}
0 & 0 & 0 \\
0 & 0 & 1\\
0 & 0 & 0
\end{array}
\right)\,. \label{notationyplusmmnus}
\end{equation}
As a direct consequence of the moving frame (\ref{framecp2}) we get
\begin{equation}
(\phi^{\dagger} y_-\phi)^{\dagger}=\phi^{\dagger} y_+\phi\,.
\end{equation}
Note that, over the space $\field{R}$, the set $\{y_-, y_+\}$ 
spans the same space as $\{s_1, s_2\}$.

Requiring that the parameterization of a surface $\mathcal{F}$ be conformal 
leads to the following conditions:
\begin{eqnarray}
&&g_{\xi\xi}=(\partial X, \partial X)=-\frac{1}{2}e^u\rm{tr}
(y_-)^2 =0\,, \nonumber \\
&&g_{\bar{\xi}\bar{\xi}}=(\bar{\partial} X, \bar{\partial} X)=
-\frac{1}{2}e^u\rm{tr}(y_+)^2 =0\,, \nonumber \\
&&g_{\xi\bar{\xi}}=(\partial X, \bar{\partial} X)=
\frac{1}{2}e^u\rm{tr}(y_-y_+)=\frac{1}{2}e^u\,,
\label{movmetricforconfor}
\end{eqnarray}
and 
\begin{eqnarray}
(\partial X, \eta_j) &=& -\frac{1}{2}e^{\frac{u}{2}}\rm{tr}(y_-s_j)=0\,,
\nonumber \\
(\bar{\partial} X, \eta_j) &=& -\frac{1}{2}e^{\frac{u}{2}}\rm{tr}(y_+s_j)=0\,,
\nonumber \\
(\eta_j, \eta_k) &=& -\frac{1}{2}\rm{tr}(s_js_k)=\delta_{jk}\,,
\label{orthonoraftermovmetricforconfor}
\end{eqnarray}
where $j, k=3, \ldots, 8$. 
Thus, we have the following proposition.
\begin{pre}
In the adjoint $SU(3)$ representation, the moving frame 
(\ref{framecp2}) of a conformally parametrized surface
$\mathcal{F}$ is described in terms of holomorphic solutions
$\{w_1, w_2\}$ of the $\field{C}P^2$ equations (\ref{cp2}) 
by the formulae
\begin{eqnarray}
\!\!\!\!&\eta_1=-\frac{i}{{A_2}^2}\left(\begin{array}{ccc}
\delta & \beta & \gamma \\
\bar{w}_1\delta & \bar{w}_1\beta & \bar{w}_1\gamma\\
\bar{w}_2\delta & \bar{w}_2\beta & \bar{w}_2\gamma
\end{array}
\right)\,,
\quad
\eta_2=-\frac{i}{{A_2}^2}\left(\begin{array}{ccc}
\bar{\delta} & w_1\bar{\delta} & w_2\bar{\delta} \\
\bar{\beta} & w_1\bar{\beta} & w_2\bar{\beta} \\
\bar{\gamma} & w_1\bar{\gamma} & w_2\bar{\gamma} \\
\end{array}
\right)\,,
\label{comparingsys1}
\end{eqnarray}
and
\begin{equation}
u=\ln(\frac{\rho}{{A_2}^{2}})\,, \label{functionu}
\end{equation}
where we define
\begin{eqnarray}
&&\delta=\bar{w}_1\partial w_1+\bar{w}_2\partial w_2\,, \nonumber \\
&&\beta=w_1\bar{w}_2 \partial w_2 - (1+|w_2|^2) \partial w_1\,,
\nonumber \\
&&\gamma=\bar{w}_1w_2 \partial w_1 - (1+|w_1|^2) \partial w_2\,,
\nonumber \\
&&\rho=|\partial w_1|^2+|\partial w_2|^2+|w_2\partial w_1-
w_1\partial w_2|^2\,.
\end{eqnarray}
\end{pre}
\vskip 5pt
\noindent
\textbf{Proof} Using the polar decomposition of the $SU(3)$ group
given by (\ref{generalsu3}), and calculating the products in the frame 
(\ref{framecp2}), yields
\begin{eqnarray}
&&\eta_1=ie^{\frac{u}{2}}\left(\begin{array}{ccc}
-a_1b_1\sin^2\alpha&-b_1\sin\alpha\,\, \zeta & -b_1\sin\alpha\,\, \mu \\
\chi\,\, a_1\sin\alpha & \chi\,\, \zeta &\chi\,\, \mu  \\
\nu\,\, a_1\sin\alpha & \nu\,\, \zeta & \nu\,\, \mu
\end{array}
\right)\,,
\nonumber \\
&&\eta_2=ie^{\frac{u}{2}}\left(\begin{array}{ccc}
-\bar{a}_1\bar{b}_1\sin^2\alpha&\bar{\chi}\,\,\bar{a}_1\sin\alpha  & 
\bar{\nu}\,\,\bar{a}_1\sin\alpha \\
-\bar{b}_1\sin\alpha\,\,\bar{\zeta} & \bar{\chi}\,\,\bar{\zeta} & 
\bar{\nu}\,\,\bar{\zeta}\\
-\bar{b}_1\sin\alpha\,\,\bar{\mu}& \bar{\chi}\,\,\bar{\mu} & 
\bar{\nu}\,\,\bar{\mu}
\end{array}
\right)\,,
\label{comparingsys2}
\end{eqnarray}
where
\begin{eqnarray}
&&\chi=-a_1b_2-\bar{a}_2b_1e^{i\varphi}\cos\alpha\,, \qquad
\zeta=-b_1\bar{b}_2+a_1a_2e^{-i\varphi}\cos\alpha\,, \nonumber \\
&&\mu=\bar{a}_2b_1+a_1b_2e^{-i\varphi}\cos\alpha\,, \qquad
\nu=a_1a_2-b_1\bar{b}_2e^{i\varphi}\cos\alpha\,.
\end{eqnarray}
Comparing (\ref{comparingsys1}) with (\ref{comparingsys2}) we obtain 
an underdetermined system of eight equations for nine unknown functions 
$a_i$, $b_i$ $\in$ $\field{C}$, $i=1,2$ and $\alpha$,
$\varphi$, $u$ $\in$ $\field{R}$. 
This system has a unique solution up to a $U(1)$ transformation. 
In other words, the phase $e^{i\varphi}$ remains arbitrary.

A straightforward algebraic computation gives $a_i$, $b_i$ and $\alpha$ 
in terms of the fields $w_1$ and $w_2$ for the $\field{C}P^2$ model.
Explicitly, 
\begin{eqnarray}
&&a_1=\frac{\sqrt{\delta\, \kappa}}{A_2\sin\alpha}e^{-{u}/{4}}\,,
\qquad
b_1=\frac{\sqrt{\delta / \kappa}}{A_2\sin\alpha}
e^{-{u}/{4}}\,,\nonumber \\
&&a_2=-\frac{e^{i\varphi}\bar{\partial} \bar{w}_2
(w_2\partial w_1-w_1\partial w_2)}{\rho\sin\alpha\cos\alpha}\,,
\quad
b_2=\frac{e^{i\varphi}\bar{\partial} \bar{w}_1
(w_2\partial w_1-w_1\partial w_2)}{\rho\sin\alpha\cos\alpha}\,,
\nonumber \\
&&\sin^2\alpha=\frac{|\partial w_1|^2+|\partial w_2|^2}{\rho}\,,
\qquad
\cos^2\alpha=\frac{|w_2\partial w_1 -w_1\partial w_2|^2}{\rho}
\,, \label{unkonwnfuncdet}
\end{eqnarray}
with $u$ as in (\ref{functionu}) and 
\begin{equation}
\kappa=\frac{\delta \cos\alpha}{w_2\partial w_1-w_1\partial w_2} 
e^{-i\varphi}\,.
\end{equation}
With the above, we can determine the moving frame (\ref{framecp2}) 
on $\mathcal{F}$, expressed in terms of the $w_1$ and $w_2,$ 
in the required form (\ref{comparingsys1}). 
That ends the proof since 
by direct computation one can check that the compatibility 
conditions, {i.e.,}
$\partial \bar{\partial}X=\bar{\partial} \partial X,$ for (\ref{framecp2}) 
are trivially satisfied. $\;\;\;\;\;\;\;\;\;\;\hspace{4cm}\square$
\vskip 4pt
\noindent
\textbf{Remark:} The explicit expressions 
for the complex normals $\eta_3, \ldots, \eta_8$ to this surface 
immersed in $su(3)$ have been calculated. 
However, the resulting expressions (in terms of $w_1$ and $w_2)$ 
are rather involved. 
A specific example is given in Appendix A.

The real-valued function $u$ given by (\ref{unkonwnfuncdet}) 
represents the total energy \cite{Zakrzewski} of the $\field{C}P^2$ model 
defined over $S^2$, since 
\begin{equation}
u=2\ln(|Dz|^2+|\bar{D}z|^2)\,,
\end{equation}
holds. 

Using the components of the induced metric (\ref{firstfundform}),  
we can write the nonzero Christoffel symbols of the second kind as
\begin{equation}
\Gamma_{11}^1=\frac{1}{q}\partial q\,, \qquad
\Gamma_{22}^2=\frac{1}{q}\bar{\partial} q\,.
\end{equation}
In this case, $q$ defined in (\ref{definedq}), becomes 
\begin{equation}
q=\frac{|\partial w_1|^2+|\partial w_2|^2+|w_1\partial w_2 -w_2
\partial w_1|^2}{2\, (1+|w_1|^2+|w_2|^2)^2}\,. \label{qforholomorcp2}
\end{equation}

Finally, taking into account (\ref{normaconditionswitdelta}), 
(\ref{movmetricforconfor}) and (\ref{orthonoraftermovmetricforconfor}), 
the moving frame (\ref{movframe}) on $\mathcal{F}$ satisfies the 
following Gauss-Weingarten equations
\begin{eqnarray}
&&\partial^2X=\frac{\partial q}{q} \partial X + J_j\eta_j\,,
\nonumber \\
&&\partial \bar{\partial} X=  H_j\eta_j\,,
\nonumber \\
&&\partial \eta_j = -2\frac{{A_2}^{\!\!2}}{\rho}(H_j \partial X
+J_j \bar{\partial} X) + S_{jk} \eta_k\,,
\label{gaussweingarten1}
\end{eqnarray}
and 
\begin{eqnarray}
&&\bar{\partial}^2X=\frac{\bar{\partial} q}{q} \bar{\partial} X 
+ \bar{J}_j\eta_j\,,
\nonumber \\
&&\bar{\partial} \partial  X=  H_j\eta_j\,,
\nonumber \\
&&\bar{\partial} \eta_j = -2\frac{{A_2}^{\!\!2}}{\rho}(\bar{J}_j 
{\partial X}+H_j \bar{\partial} X) + \bar{S}_{jk} \eta_k\,,
\label{gaussweingarten2}
\end{eqnarray}
where
\begin{equation}
J_j=-\frac{1}{2}\rm{tr}(\partial^2X\,\eta_j)\,, \qquad
H_j=-\frac{1}{2}\rm{tr}(\partial \bar{\partial} X\,\eta_j)\,,
\label{defgaussweingarten1}
\end{equation}
and 
\begin{equation}
S_{jk}+S_{kj}=0\,, \qquad \bar{S}_{jk}+\bar{S}_{kj}=0\,,
\qquad j \ne k = 3, \ldots, 8\,.
\label{defgaussweingarten2}
\end{equation}

The Gauss-Codazzi-Ricci equations, which are the compatibility conditions 
for (\ref{gaussweingarten1}) and (\ref{gaussweingarten2}), coincide with 
the equations of the $\field{C}P^{N-1}$ model. 
However, the explicit forms of the coefficients $H_j$ and $J_j$ depend 
locally on the chosen orthonormal basis $\{\eta_3,\ldots,\eta_8\}$ 
of the space normal to the surface $\mathcal{F}$ at a given point 
$p=(\xi_0, \bar{\xi}_0)\in X$. Note that quantities $H_j$ and $J_j$ 
are not completely arbitrary. 
Using (\ref{movmetricforconfor}) and the fact that $J=0$, 
it becomes clear that the complex tangent vectors have to satisfy 
the following differential constraints
\begin{equation}
(\partial^2X\,, \bar{\partial} \partial  X)=0\,,
\qquad
(\bar{\partial}^2X\,, \partial \bar{\partial} X)=0\,.
\end{equation}

For any holomorphic solution $(w_i,\bar{w}_i)$ $i=1,2$ of the 
$\field{C}P^2$ model, we computed explicitly the form of the 
first and second fundamental forms, $I$ and $II,$ and the mean curvature 
vector $\mathcal{H}$ of a conformally parametrized surface $\mathcal{F}$ 
at some regular point $p=(\xi_0, \bar{\xi}_0)\in X.$
They are
\begin{eqnarray}
&&I=\frac{\rho}{{A_2}^{\!\!2}} d\xi d\bar{\xi}\,, \nonumber \\
&&II=\left(\partial^2 X-\frac{\partial q}{q}\,\partial X\right)d\xi^2+2
\partial\bar{\partial}Xd\xi d\bar{\xi}+
\left(\bar{\partial}^2 X-\frac{\bar{\partial} q}{q}\,\bar{\partial} X\right) 
d\bar{\xi}^2\,,
\nonumber \\
&&\mathcal{H}=\frac{2}{q}\partial\bar{\partial}X\,,
\end{eqnarray}
respectively. 
The second derivatives of the Weierstrass representation 
$X$ can be computed from (\ref{comparingsys2}).

One can also compute some of the global properties of surfaces associated 
with the $\field{C}P^2$ sigma model, using the well-known formulae 
(see {e.g.,} \cite{KobayashiNomizu, Willmore}). 
For instance, for any set of holomorphic solutions 
($w_i, \bar{w}_i$), $i=1,2$, of the $\field{C}P^2$ model, the Willmore 
functional assumes the form
\begin{equation}
W=-4i\int_{\Omega}\frac{1}{q}[\partial P, \bar{\partial} P]^2
d\xi d\bar{\xi}\,,
\end{equation}
whose values depend only on the fields and their derivatives 
on the boundary $\partial \Omega$ of the open set $\Omega$.

Under the above assumptions and provided that the $\field{C}P^2$ model is 
defined on the whole Riemannian sphere $S^2$, the topological charge becomes
\begin{equation}
Q=-\frac{1}{8\pi}\int_{S^2} q \, d\xi d\bar{\xi}\,.
\label{topologicalcahrge}
\end{equation}
If the above integral exists, then it is an integer which globally 
characterizes the surface.

\section{The Weierstrass formula for immersion of surfaces in the 
$su(2)$ and $su(3)$ algebras \label{Weierstrassformulaforsu1su2}}
In this section we apply the general idea of Weierstrass representation 
of surfaces given in Section \ref{secweierstrass} to two specific cases, 
namely, the $\field{C}P^1$ and $\field{C}P^2$ models. 
For each case, we first find the concrete form of the generalized 
Weierstrass representation of surfaces associated with these models 
and then we give the corresponding Weierstrass data for the 
holomorphic solutions. 

It is known \cite{Grundlandetal,GrundlandSnobl} that, with the projector 
$P$ given by (\ref{projector}), one can compute explicitly the formula 
for immersion (\ref{intforimm}) in terms of the complex fields $w_i$ 
of the equations of motion of the model.

We start with the case $N=2$. The orthogonal projector $P$ and matrix 
$K$ are then given by 
\begin{equation}
P=I_2-\frac{1}{A_1} \left(\begin{array}{cc}
1 & w_1  \\
\bar{w}_1 & w_1\bar{w}_1  \\
\end{array}
\right)\,, \label{projectorcp1}
\end{equation}
and
\begin{equation}
K=\frac{1}{{A_1}^{\!\!2}} \left(\begin{array}{cc}
\bar{w}_1 \bar{\partial} w_1-w_1\bar{\partial} \bar{w}_1  & -(\bar{\partial}
 w_1+w_1^2 \bar{\partial} \bar{w}_1)  \\
(\bar{\partial} \bar{w}_1+\bar{w}_1^2 \bar{\partial} {w}_1) & w_1
\bar{\partial} \bar{w}_1 -\bar{w}_1 \bar{\partial} w_1  \\
\end{array}
\right)\,, \label{matrixforcp1}
\end{equation}
where, as usual, $w_1$ is the homogeneous coordinate defined 
by (\ref{coordinw}). 
Based on the expression of the matrix $K$ for the $\field{C}P^1$ model, 
the Weierstrass data follows from (\ref{complex1formimmer}). 
In order to obtain real-valued $1$-forms we decompose $dX$ given in 
(\ref{complex1formimmer}) into its real and imaginary parts, 
\begin{equation}
dX=dX^1+idX^2\,.
\end{equation}
So,
\begin{eqnarray}
dX^1=\frac{i}{2}\Big[(K^{\dagger}-\bar{K})d\xi + (K-K^T)d\bar{\xi}\Big]
\,, \nonumber \\
dX^2=\frac{1}{2}\Big[(K^{\dagger}+\bar{K})d\xi + (K+K^T)d\bar{\xi}\Big]
\,. \label{real1formscp1}
\end{eqnarray}
It is easily seen that $dX^1$ is skew-symmetric and $dX^2$ is symmetric. 
Realizing that the $2D$ surface associated with the $\field{C}P^1$ model 
is immersed in the $su(2)$ algebra, the two real-valued $1$-forms can 
also be expressed in terms of the Pauli matrices. 
Since $dX^1$ is skew-symmetric and $dX^2$ is symmetric, the $1$-forms 
can be represented as
\begin{equation}
dX^1=idX_2\sigma_2\,, \qquad dX^2=dX_1\sigma_1+dX_3\sigma_3\,,
\label{decompositioncp1}
\end{equation}
where $\sigma_1$, $\sigma_2$ and $\sigma_3$ are the Pauli matrices
\begin{equation}
\sigma_1=\left(\begin{array}{cc}
0 & 1  \\
1 & 0  \\
\end{array}
\right)\,, \qquad
\sigma_2=\left(\begin{array}{cc}
0 & -i  \\
i & 0  \\
\end{array}
\right)\,, \qquad
\sigma_3=\left(\begin{array}{cc}
1 & 0  \\
0 & -1  \\
\end{array}
\right)\,. \label{sigmamatrices}
\end{equation}
After substituting the matrix $K$ from (\ref{matrixforcp1}) into  
(\ref{real1formscp1}) and comparing the result with 
(\ref{decompositioncp1}), it is easy to see that the real-valued 
$1$-forms $dX_i$, $i=1,2,3$, can be expressed in terms of the solutions 
of the Euler-Lagrange equations of the $\field{C}P^1$ model. 
Indeed,
\begin{eqnarray}
&&dX_1=\frac{1}{2{A_1}^{\!\!2}} \left(\Big[(1-\bar{w}_1^2)\partial w_1 
-(1-w_1^2) \partial \bar{w}_1\Big]d\xi + \rm{c.c.} \right)\,, 
\nonumber \\
&&dX_2=\frac{i}{2{A_1}^{\!\!2}} \left(\Big[(1+{w}_1^2)\partial \bar{w}_1 
+(1+\bar{w}_1^2) \partial {w}_1\Big]d\xi - \rm{c.c.} \right)\,, 
\nonumber \\
&&dX_3=\frac{1}{{A_1}^{\!\!2}} \left(\Big[w_1\partial \bar{w}_1-\bar{w}_1 
\partial w_1 \Big]d\xi + \rm{c.c.} \right)\,, 
\label{weiestrassforcp1}
\end{eqnarray}
where ``c.c." denotes the complex conjugate.
In fact, these real-valued $1$-forms constitute the generalized 
Weierstrass formula for immersion for the $\field{C}P^1$ model. 

Now, we further restrict ourselves to the holomorphic solutions of the 
$\field{C}P^1$ model. 
This restriction is necessary if the model is defined on $S^2$ with a 
finite action \cite{Zakrzewski}.
Using holomorphic solutions, $dX_i,$ $i=1,2,3$, can be reduced into
\begin{eqnarray}
&&dX_1=\frac{1}{2}\, \partial \left(\frac{w_1+\bar{w}_1}{A_1}\right) 
d\xi +\rm{c.c.} \,, 
\nonumber \\
&&dX_2=\frac{i}{2}\left[\partial \left(\frac{w_1-\bar{w}_1}{A_1}\right) 
d\xi -\rm{c.c.} \right] \,, 
\nonumber \\
&&dX_3=-\partial \left(\frac{|w_1|^2}{A_1}\right) d\xi +\rm{c.c.}\,. 
\label{weiestrassforcp1forholomorcp1}
\end{eqnarray}
Integration gives
\begin{equation}
X_1=\frac{w_1+\bar{w}_1}{2\,A_1}\,, \qquad 
X_2=i\frac{w_1-\bar{w}_1}{2\,A_1}\,, \qquad
X_3=-\frac{|w_1|^2}{A_1}\,,
\label{coordinatesofrcp1}
\end{equation}
where the constants of integration are set to zero.

It is well-known that the $2D$ surface associated with the holomorphic 
solutions of the $\field{C}P^1$ model is the surface of a 
sphere \cite{Zakrzewski}. 
Confirmation of that result follows from (\ref{coordinatesofrcp1}). 
Indeed, upon elimination of $w_1$ and ${\bar{w}_1},$ we obtain 
\begin{equation}
X_1^2+X_2^2+\left(X_3+\frac{1}{2}\right)^2=\frac{1}{4}\,.
\label{eqsurfaceforcp1}
\end{equation}
So, all points of the $2D$ surface lie on the surface of a sphere of 
radius $1/2$ centered at $(0,0,-1/2)$.

We now consider the case $N=3$. The corresponding orthogonal projector 
$P$ is given in (\ref{projectorcp2}) and matrix $K = -i\eta_2$ with 
$\eta_2$ in (\ref{comparingsys1}). 
Since the $2D$ surface associated with the $\field{C}P^2$ model is 
immersed in the $su(3)$ algebra, the two real-valued $1$-forms, $dX^1$ and 
$dX^2$, obtained by decomposing $dX=i(K^{\dagger}d\xi + K d\bar{\xi})$ 
into real and imaginary parts, can be expressed in terms of the orthonormal 
basis of the Lie algebra $su(3)$. 
Keeping in mind that $dX^1$ is skew-symmetric and $dX^2$ is symmetric, 
the real-valued $1$-forms are given by
\begin{eqnarray}
&&dX^1=dX_2 s_2 + dX_5 s_5 + dX_6 s_6\,, \nonumber \\
&&dX^2=i \big(dX_1 s_1 + dX_3 s_3 + dX_4 s_4 + dX_7 s_7 + dX_8 s_8 \big)\,,  
\label{decompositioncp2}
\end{eqnarray}
where the Gell-Mann matrices $s_i$, $i=1,\ldots, 8$, are given in 
(\ref{matricess8}).

Using $K=-i\eta_2$ and comparing (\ref{real1formscp1}) with 
(\ref{decompositioncp2}), it follows that 
the real-valued $1$-forms $dX_i, i=1,\dots, 8$, can be expressed 
in terms of the solutions of the Euler-Lagrange equations of the 
$\field{C}P^2$ model as 
\begin{eqnarray}
dX_1&=&\frac{1}{2{A_2}^{\!\!2}}\Big(\big[(w_2^2-w_1^2)(\bar{w}_1\partial 
\bar{w}_2-\bar{w}_2\partial \bar{w}_1)-(\bar{w}_2^2-\bar{w}_1^2)({w}_1
\partial {w}_2-{w}_2\partial {w}_1) \nonumber \\
&\,&-w_2\partial \bar{w}_1+\bar{w}_2\partial {w}_1-w_1\partial \bar{w}_2
+\bar{w}_1\partial {w}_2\big] d\xi + \rm{c.c.} \Big)\,,
\nonumber \\
dX_2&=&\frac{i}{2{A_2}^{\!\!2}}\Big(\big[(w_1^2+w_2^2)(\bar{w}_2\partial 
\bar{w}_1-\bar{w}_1\partial \bar{w}_2)+(\bar{w}_1^2+\bar{w}_2^2)({w}_2
\partial {w}_1-{w}_1\partial {w}_2) \nonumber \\
&\,&+w_2\partial \bar{w}_1+\bar{w}_2\partial {w}_1-w_1\partial \bar{w}_2-
\bar{w}_1\partial {w}_2\big] d\xi - \rm{c.c.} \Big)\,,
\nonumber \\
dX_3&=&\frac{1}{2{A_2}^{\!\!2}}\Big(\big[w_2\partial \bar{w}_2-{w}_1\partial 
\bar{w}_1-\bar{w}_2\partial {w}_2+\bar{w}_1\partial {w}_1 
\nonumber \\
&\,&+ 2 |w_1|^2 ({w}_2\partial \bar{w}_2-\bar{w}_2\partial {w}_2) 
- 2 |w_2|^2 ({w}_1\partial \bar{w}_1-\bar{w}_1\partial {w}_1) \big] d\xi + 
\rm{c.c.} \Big)\,,
\nonumber \\
dX_4&=&\frac{\sqrt{3}}{2{A_2}^{\!\!2}}\Big(\big[w_1\partial \bar{w}_1+{w}_2
\partial \bar{w}_2-\bar{w}_1\partial {w}_1-\bar{w}_2\partial {w}_2 \big] 
d\xi + \rm{c.c.} \Big)\,,
\nonumber \\
dX_5&=&-\frac{i}{2{A_2}^{\!\!2}}\Big(\big[(1+\bar{w}_1^2+|w_2|^2)
\partial w_1 + (1+{w}_1^2+|w_2|^2) \partial \bar{w}_1
\nonumber \\
&\,&+({w}_2\partial \bar{w}_2-\bar{w}_2\partial {w}_2)(w_1-\bar{w}_1)\big] 
d\xi - \rm{c.c.} \Big)\,,
\nonumber \\
dX_6&=&-\frac{i}{2{A_2}^{\!\!2}}\Big(\big[(1+\bar{w}_2^2+|w_1|^2)
\partial w_2 + (1+{w}_2^2+|w_1|^2) \partial \bar{w}_2
\nonumber \\
&\,&+({w}_1\partial \bar{w}_1-\bar{w}_1\partial {w}_1)(w_2-\bar{w}_2)\big] 
d\xi - \rm{c.c.} \Big)\,,
\nonumber \\
dX_7&=&\frac{1}{2{A_2}^{\!\!2}}\Big(\big[(1-{w}_1^2+|w_2|^2)\partial 
\bar{w}_1 - (1-\bar{w}_1^2+|w_2|^2) \partial {w}_1
\nonumber \\
&\,&+(\bar{w}_2\partial {w}_2-{w}_2\partial \bar{w}_2)(w_1+\bar{w}_1)\big] 
d\xi + \rm{c.c.} \Big)\,,
\nonumber \\
dX_8&=&\frac{1}{2{A_2}^{\!\!2}}\Big(\big[(1-{w}_2^2+|w_1|^2)\partial 
\bar{w}_2 - (1-\bar{w}_2^2+|w_1|^2) \partial {w}_2
\nonumber \\
&\,&+(\bar{w}_1\partial {w}_1-{w}_1\partial \bar{w}_1)(w_2+\bar{w}_2)\big] 
d\xi + \rm{c.c.} \Big)\,. 
\label{weiestrassforcp2}
\end{eqnarray}
These eight real-valued $1$-forms constitute the generalized 
Weierstrass formula for immersion for the $\field{C}P^2$ model. 
\vskip 4pt
\noindent
\textbf{Remark:}
Note that the reflection transformations in independent or dependent 
variables and their complex conjugates preserve the form of the 
$\field{C}P^2$ model. 
So does the generalized $SU(2)$ transformation.
Indeed, if the complex-valued functions $u_1$ and $u_2$ are solutions of 
the $\field{C}P^2$ model, then the complex-valued functions 
$w_1$ and $w_2$ defined by the generalized $SU(2)$ transformation,
\begin{eqnarray}
w_1\rightarrow \frac{a^2u_1-b^2u_2-\sqrt{2}\,a\,b}{\sqrt{2}(a\,
\bar{b}\,u_1+\bar{a}\,b\,u_2)+|a|^2-|b|^2}\,, \nonumber \\
w_2\rightarrow \frac{-\bar{b}^2u_1+\bar{a}^2u_2-\sqrt{2}\,\bar{a}\,
\bar{b}}{\sqrt{2}(a\,\bar{b}\,u_1+\bar{a}\,b\,u_2)+|a|^2-|b|^2}\,,
\end{eqnarray}
for $a, b \in \field{C}$ such that $|a|^2+|b|^2=1$, are also solutions 
of the $\field{C}P^2$ model.

These transformations can be used to restrict the range of parameters 
appearing in the explicit form of solutions of the $\field{C}P^2$ model. 
They allow one to simplify the Weierstrass representation.

Again, we restrict ourselves to the holomorphic solutions 
of the $\field{C}P^2$ model.
In that case, the eight real-valued $1$-forms $dX_i, i=1,\ldots, 8$, are
\begin{eqnarray}
&&dX_1=\frac{1}{2}\,\partial \left(\frac{w_1\bar{w}_2+\bar{w}_1 w_2}{A_2}
\right)d \xi + \rm{c.c.}\,,  \nonumber \\
&&dX_2=\frac{i}{2}\left[\partial \left(\frac{w_1\bar{w}_2-\bar{w}_1 w_2}{A_2}
\right)d \xi - \rm{c.c.} \right]\,,  \nonumber \\
&&dX_3=\frac{1}{2}\,\partial \left(\frac{|w_1|^2-|w_2|^2}{A_2}\right)d \xi + 
\rm{c.c.}\,, \nonumber \\
&&dX_4=-\frac{\sqrt{3}}{2}\,\partial \left(\frac{|w_1|^2+|w_2|^2}{A_2}\right)d 
\xi + \rm{c.c.}\,, \nonumber \\
&&dX_5=-\frac{i}{2}\left[\partial \left(\frac{w_1-\bar{w}_1}{A_2}\right)d \xi- 
\rm{c.c.} \right]\,,  \nonumber \\
&&dX_6=-\frac{i}{2}\left[\partial \left(\frac{w_2-\bar{w}_2}{A_2}\right)d \xi- 
\rm{c.c.} \right]\,,  \nonumber \\
&&dX_7=-\frac{1}{2}\,\partial \left(\frac{w_1+\bar{w}_1}{A_2}\right)d \xi + 
\rm{c.c.} \,,  \nonumber \\
&&dX_8=-\frac{1}{2}\,\partial \left(\frac{w_2+\bar{w}_2}{A_2}\right)d \xi + 
\rm{c.c.} \,.  
\label{weiestrassforcp1forholomorcp2}
\end{eqnarray}
Ignoring integration constants, after integration we obtain
\begin{eqnarray}
&&X_1=\frac{w_1\bar{w}_2+\bar{w}_1w_2}{2\,A_2}\,, \qquad 
X_2=i\frac{w_1\bar{w}_2-\bar{w}_1w_2}{2\,A_2}\,, \qquad 
X_3=\frac{|w_1|^2-|w_2|^2}{2\,A_2}\,, \nonumber \\
&&X_4=-\sqrt{3}\frac{|w_1|^2+|w_2|^2}{2\,A_2}\,, \qquad 
X_5=-i\frac{w_1-\bar{w}_1}{2\,A_2}\,, \qquad
X_6=-i\frac{w_2-\bar{w}_2}{2\,A_2}\,, \nonumber \\
&&X_7=-\frac{w_1+\bar{w}_1}{2\,A_2}\,, \qquad
X_8=-\frac{w_2+\bar{w}_2}{2\,A_2}\,,
\label{coordinatesofrcp2}
\end{eqnarray}
which determines the coordinates of the radius vector 
$\vec{X}=(X_1,\dots, X_8)$ of a two-dimensional surface in $\field{R}^8$.

Note that in the limiting cases $w_i\rightarrow w/\sqrt{2}$, $i=1,2$, 
or $w_1\rightarrow 0$, or $w_2\rightarrow 0$, the generalized 
Weierstrass formula (\ref{weiestrassforcp2}) for immersion of 
the $\field{C}P^2$ model reduces (after straightforward manipulations)
to the generalized Weierstrass formula (\ref{weiestrassforcp1}) 
for immersion of the $\field{C}P^1$ model.  
Consequently, the coordinates of radius vector $\vec{X}$ in 
(\ref{coordinatesofrcp2}) for the holomorphic solutions of the $\field{C}P^2$ 
model then reduce to the coordinates of $\vec{X}$ in (\ref{coordinatesofrcp1}) 
for the holomorphic solutions of the $\field{C}P^1$ model.

When dealing with the $2D$ surface associated with the holomorphic 
solutions of the $\field{C}P^2$ model, all points lie on the affine sphere, 
\begin{equation}
4X_1^2+4X_2^2+4X_3^2+\frac{2}{\sqrt{3}}X_4+X_5^2+X_6^2+X_7^2+X_8^2=0\,.
\label{surfaceeqforcp2}
\end{equation}
It is straightforward to show that the coordinates of the radius vector
(\ref{coordinatesofrcp2}) satisfy (\ref{surfaceeqforcp2}).

\section{Examples of surfaces associated with the $\field{C}P^{N-1}$ sigma models 
\label{secexamples}}
Using elementary examples, we will illustrate the concept of constructing 
surfaces associated with the $\field{C}P^{N-1}$ model.

\subsection{Examples of holomorphic solutions of the $\field{C}P^2$ sigma model} 
From the form of the $\field{C}P^2$ model, it is readily seen that the 
holomorphic functions are solutions of the $\field{C}P^2$ model. 
We now concentrate on the following class of holomorphic solutions of the 
$\field{C}P^2$ model:
\begin{equation}
w_1=a_1\xi^m\,, \qquad w_2=a_2\xi^n\,, \label{holsolutionscp2examp}
\end{equation}
where $a_1$ and $a_2$ are complex constants and $m$ and $n$ are real 
constants. 
For holomorphic solutions $J=0$ and the induced metric is conformal.
Using the solutions in (\ref{holsolutionscp2examp}), that metric is given by 
\begin{equation}
I=\frac{|a_1|^2 |\xi|^{2m}(m^2+|a_2|^2(m-n)^2 |\xi|^{2n})+|a_2|^2 n^2 
|\xi|^{2n}}{|\xi|^2(1+|a_1|^2|\xi|^{2m}+|a_2|^2|\xi|^{2n})^2}d\xi d\bar{\xi} 
\,. \label{conmetexamp}
\end{equation}
The Gaussian curvature $\mathcal{K}$ is computed from (\ref{gaussiancur2}).
After simplification, 
\begin{equation}
\mathcal{K}=4-\frac{2|a_1|^2 |a_2|^2 m^2 n^2 (m-n)^2 |\xi|^{2(m+n)} 
(1+|a_1|^2|\xi|^{2m}+|a_2|^2|\xi|^{2n})^3}{\Big(|a_1|^2 |\xi|^{2m}
(m^2+|a_2|^2(m-n)^2 |\xi|^{2n})+|a_2|^2 n^2 
|\xi|^{2n}\Big)^3}\,. \label{Gaussiancurexampgeneral}
\end{equation}
In general, $\mathcal{K}$ is not constant. 
However, $\mathcal{K}$ is constant for certain values of 
$a_1$, $a_2$, $m$ and $n$. 
For example, if the second term in (\ref{Gaussiancurexampgeneral}) 
vanishes or equals to a constant, 
then the surfaces associated with the holomorphic solutions 
(\ref{holsolutionscp2examp}) of the $\field{C}P^2$ model will have 
constant Gaussian curvature. 
This happens when 
\begin{enumerate}
\item[(i)] $a_1=0$, $a_2=0$, $m=0$, $n=0$ and $m=n$ or a combination thereof.
For these choices the second term in (\ref{Gaussiancurexampgeneral}) vanishes; 
or
\item[(ii)] $n=2m$ and $|a_1|^2=\pm 2|a_2|$ simultaneously. 
The second term in (\ref{Gaussiancurexampgeneral}) then reduces to a constant.
\end{enumerate}

Not surprisingly, constant Gaussian curvature occurs when $a_1=0$ or $a_2=0$ 
because the $\field{C}P^2$ model then reduces to the $\field{C}P^1$ model.
Hence, the surfaces must have constant Gaussian curvature. 

We now consider a case of constant Gaussian curvature surfaces 
associated with specific holomorphic solutions (\ref{holsolutionscp2examp}) 
of the $\field{C}P^2$ model. 
For simplicity, we take
\begin{equation}
w_1=\xi\,, \qquad w_2=\frac{1}{2}\xi^2\,. \label{holsolutionscp2special}
\end{equation}
The first fundamental form and the Gaussian curvature then are
\begin{eqnarray}
&&I=\frac{4}{(2+|\xi|^2)^2}d\xi d\bar{\xi}\,, \nonumber \\
&&\mathcal{K}=2\,. \label{expspecialfirstGauss}
\end{eqnarray}
Upon substitution of (\ref{holsolutionscp2special}) into 
(\ref{coordinatesofrcp2}), the coordinates of the radius vector $\vec{X}$ 
become 
\begin{eqnarray}
&&X_1=\frac{|\xi|^2(\xi+\bar{\xi})}{(2+|\xi|^2)^2}\,,\quad
X_2=i\frac{|\xi|^2(\bar{\xi}-{\xi})}{(2+|\xi|^2)^2}\,,\quad
X_3=\frac{|\xi|^2(4-|\xi|^2)}{2\,(2+|\xi|^2)^2}\,, \nonumber \\
&&X_4=-\frac{\sqrt{3}}{2}\left(1-\frac{4}{(2+|\xi|^2)^2}\right)\,, \qquad
X_5=-i\frac{2(\xi-\bar{\xi})}{(2+|\xi|^2)^2}\,, \label{coorradiiusexamp} \\
&&X_6=-i\frac{(\xi^2-\bar{\xi}^2)}{(2+|\xi|^2)^2}\,, \quad
X_7=-\frac{2(\xi+\bar{\xi})}{(2+|\xi|^2)^2}\,, \quad
X_8=-\frac{(\xi^2+\bar{\xi}^2)}{(2+|\xi|^2)^2}\,. \nonumber
\end{eqnarray}
Of course, the above coordinates satisfy the relation (\ref{surfaceeqforcp2}). 
Hence, the surface associated with the specific solutions 
(\ref{holsolutionscp2special}) of the $\field{C}P^2$ model is an affine 
sphere. 

\subsection{Mixed solutions of the $\field{C}P^2$ sigma model} 
In this subsection we analyze the mixed solutions of the 
$\field{C}P^2$ model and give the first fundamental form, 
Gaussian curvature and the Weierstrass data for a specific example. 
It is well-known \cite{Zakrzewski} that if the $\field{C}P^2$ model 
is defined over $S^2$ and the 
finiteness of the action (\ref{action2}) is required, 
then the solutions of the $\field{C}P^2$ model split into three cases: 
holomorphic solutions, anti-holomorphic solutions and mixed ones. 
Among these, the mixed solutions can be constructed either from the 
holomorphic or anti-holomorphic solutions according to the following 
procedure \cite{Grundlandetal, Zakrzewski}.

Consider three arbitrary holomorphic functions $g_i=g_i(\xi)$, $i = 1,2,3,$ 
and define the Wronskian
\begin{equation}
G_{ij}=g_i\partial g_j-g_j\partial g_i\,, \qquad i=1,2,3\,,
\end{equation}
based on any pair. 
It can be verified that the functions 
\begin{equation}
f_i=\sum_{k\ne i}^3 \bar{g}_k G_{ki}\,, \qquad i=1,2,3\,,
\end{equation}
are solutions of the $\field{C}P^2$ model.
The mixed solutions are associated with the ratios
\begin{equation}
w_1=\frac{f_1}{f_3}\,, \qquad w_2=\frac{f_2}{f_3}\,.
\end{equation}

Likewise, mixed solutions can be obtained from anti-holomorphic solutions 
by using $\bar{\partial}$ instead of $\partial$.  

We now continue with the holomorphic functions
\begin{equation}
g_1=1\,,\qquad g_2=\rm{sech}(\xi)\,, \qquad 
g_3=\tanh(\xi)\,.
\end{equation}
Using the above procedure, the mixed solutions of the $\field{C}P^2$ model 
are
\begin{equation}
w_1=\tanh(\frac{\xi-\bar{\xi}}{2})\,, \qquad 
w_2=-\frac{\tanh(\xi)+\tanh(\bar{\xi})}{\rm{sech}(\xi)+\rm{sech}(\bar{\xi})}
\,, \label{solitonsolutions}
\end{equation}
which are of soliton-type. 
These fields satisfy the equations of the $\field{C}P^2$ model.
$J=0$ for this case, as can be readily verified. 
Hence, the induced metric is conformal and given by 
\begin{equation}
I=\frac{2}{1+\cosh(\xi+\bar{\xi})}d\xi d\bar{\xi}\,.
\end{equation}

Note that holomorphicity of the solutions of the $\field{C}P^{N-1}$ model 
implies that $J=0.$ 
The converse is false as seen from the above example (\ref{solitonsolutions}). 

The Gaussian curvature is computed from the formula given in 
(\ref{gaussiancur2}) (since $J=0$) and found to be 
\begin{equation}
\mathcal{K}=1\,.
\end{equation}
After substituting the solutions (\ref{solitonsolutions}) into 
(\ref{weiestrassforcp2}) for the $\field{C}P^2$ model, 
the Weierstrass representation becomes
\begin{eqnarray}
&&dX_1=-\frac{\sinh(\bar{\xi})}{1+\cosh(\xi+\bar{\xi})}d\xi+\rm{c.c.}\,, 
\quad
dX_6=i\left[\frac{\cosh(\bar{\xi})}{1+\cosh(\xi+\bar{\xi})}d\xi-
\rm{c.c.}\right]\,, \nonumber \\
&&\;\;\;\;\;\;\;\;\;\;\;\;\;\;\;\;\;\;\;\;\;\;\;\;\;\;\;\;\;\;\;\;
dX_7=-\frac{1}{1+\cosh(\xi+\bar{\xi})}d\xi+\rm{c.c.}\,,
\label{Weierstrassforsoliton}
\end{eqnarray}
and 
\begin{equation}
dX_2=0\,,\quad dX_3=0\,,\quad dX_4=0\,,\quad dX_5=0\,,\quad dX_8=0\,.
\end{equation}
Integrating (\ref{Weierstrassforsoliton}), we obtain the coordinates of 
the radius vector $\vec{X}$:
\begin{eqnarray}
&&X_1=\rm{sech}\left(\frac{\xi+\bar{\xi}}{2}\right)\cosh\left(\frac{\xi-
\bar{\xi}}{2}\right)\,, \nonumber \\
&&X_6=i\,\rm{sech}\left(\frac{\xi+\bar{\xi}}{2}\right)\sinh\left(\frac{\xi-
\bar{\xi}}{2}\right)\,, \nonumber \\
&&X_7=-\tanh\left(\frac{\xi+\bar{\xi}}{2}\right)\,,
\end{eqnarray}
They satisfy $X_1^2+X_6^2+X_7^2=1$. 
Hence, the constant Gaussian curvature surface associated with the
soliton-like solutions (\ref{solitonsolutions}) of the $\field{C}P^2$ model 
is really immersed in $\field{R}^3$ which, in turn, corresponds to the 
immersion of the $\field{C}P^2$ model into the $\field{C}P^1$ model.

\subsection{Examples of surfaces in the $su(N)$ algebra}
We briefly discuss the non-splitting solutions $(w_i, \bar{w}_i)$, 
$i=1,\ldots, N-1$ of the $\field{C}P^{N-1}$ model invariant under the 
scaling symmetries $\{S_i\}$ as given in (\ref{cpNsymm}). 
To do so, we subject system (\ref{cpn}) to $N-1$ algebraic constraints
\begin{equation}
w_i\bar{w}_i = D_i \in \field{R}\,, \qquad i=1,\ldots, N-1\,.
\label{subjectedconstraint}
\end{equation}
If, for simplicity, we choose $D_i=1$, then the simplest solutions of this 
type are 
\begin{equation}
w_i=\frac{F_i(\xi)}{\bar{F}_i(\bar{\xi})}\,, \qquad i=1,\ldots, N-1\,,
\label{simplestsolutions}
\end{equation}
where $F_i$ and $\bar{F}_i$ are arbitrary complex-valued functions of one 
complex variable each. 
Substituting (\ref{simplestsolutions}) into (\ref{cpn}), we obtain 
a class of non-splitting  solutions of the $\field{C}P^{N-1}$ 
model which depend on one arbitrary complex-valued function of one 
variable $\xi$ and its conjugate. 
Indeed, 
\begin{equation}
w_1=\frac{F_1(\xi)}{\bar{F}_1(\bar{\xi})}\,, \qquad 
w_{j+1}=\frac{c_j}{\bar{c}_j}\frac{F_1(\xi)^{e^{i\psi}}}
{\bar{F}_1(\bar{\xi})^{e^{-i\psi}}}\,,
\qquad \qquad j=1,\ldots, N-2\,,
\label{generalsolutionsfornonsplitting}
\end{equation}
where $c_j$, $\bar{c}_j$ are complex constants and 
\begin{equation}
\psi=\pm \frac{\pi}{3}+2\pi m\,, \qquad m\in \field{Z}\,.
\label{equationforpsi}
\end{equation}
For brevity, from now on we suppress the subscript $1$ and also 
the arguments of the functions $F$ and $ {\bar{F}}.$
For this class of non-splitting solutions, 
the induced metric $g_{ij}$ has the following components 
\begin{equation}
g_{\xi\xi}=-\frac{N-3}{N^2}\frac{(F^{\prime})^2}{F^{\,2}}\,,
\quad
g_{\bar{\xi}\bar{\xi}}=-\frac{N-3}{N^2}\frac{(\bar{F}^{\prime})^2}
{\bar{F}^{\,2}}\,,
\quad
g_{\xi\bar{\xi}}=\frac{2N-3}{N^2}\frac{|F^{\prime}|^2}{|F|^{2}}\,,
\label{indmetricfornonsplitting}
\end{equation}
where prime denotes differentiation with respect to the argument.
The determinant of the induced metric then is 
\begin{equation}
g=-\frac{3(N-2)}{N^3}\frac{|F^{\prime}|^4}{|F|^{4}}\,.
\label{detformetfornonsplitting}
\end{equation}

Two interesting examples occur when $N=2$ or $N=3$. 
For $N=2$, the determinant of the  induced metric vanishes. 
Hence, the associated surface for the $\field{C}P^1$ model, subject to the 
DCs in (\ref{subjectedconstraint}), reduces to a curve in $\field{R}^3$. 
For $N=3$, the diagonal components of the induced metric vanish 
(since $J=0$). 
Hence, we have a conformal metric for the $\field{C}P^2$ model 
subject to the DCs in (\ref{subjectedconstraint}).

From (\ref{gaussiancur1}) and (\ref{gaussiancur2}) it is straightforward 
to show that the Gaussian curvature vanishes for the associated surfaces of 
the $\field{C}P^{N-1}$ model ($N\geq 3$), subject to the DCs  
(\ref{subjectedconstraint}). 
Thus, we conclude that for $N\geq 3$ the surfaces associated with 
solutions of the $\field{C}P^{N-1}$ model, which are 
invariant under dilations, always have zero Gaussian curvature, 
{i.e.,}
\begin{equation}
\mathcal{K}=0\,.
\end{equation}

Finally, let us give the coordinates of the radius vector $\vec{X}$ 
for the non-splitting solutions of the $\field{C}P^2$ model. 
After substituting the non-splitting solutions 
(\ref{generalsolutionsfornonsplitting}) of the $\field{C}P^2$ model 
into the Weierstrass representation 
(\ref{weiestrassforcp2}) and subsequent integration, the coordinates 
of the radius vector $\vec{X}$ in $\field{R}^8$ are
\begin{eqnarray}
X_1&=&\frac{i}{6\sqrt{3}\,|c|^2}|F|^{-2e^{i\psi}}(\bar{c}^2F-c^2
\bar{F}|F|^{2i\sqrt{3}})\,, \nonumber \\
X_2&=&-\frac{1}{6\sqrt{3}\,|c|^2}|F|^{-2e^{i\psi}}(\bar{c}^2F+c^2
\bar{F}|F|^{2i\sqrt{3}})\,,\nonumber \\
X_3&=&\frac{1}{6}\big((1-i\sqrt{3})\rm{ln}F+(1+i\sqrt{3})\rm{ln}
\bar{F}\big)\,,\nonumber \\
X_4&=&-\frac{1}{6}\big((i+\sqrt{3})\rm{ln}F+(-i+\sqrt{3})\rm{ln}
\bar{F}\big)\,,\nonumber \\
X_5&=&-\frac{F^2+\bar{F}^2}{6\sqrt{3}\, |F|^2}\,,\nonumber \\
X_6&=&\frac{1}{6\sqrt{3}\,|c|^2}|F|^{-2e^{i\psi}}(\bar{c}^2\bar{F}+
c^2{F}|F|^{2i\sqrt{3}})\,,\nonumber \\
X_7&=&\frac{i(F^2-\bar{F}^2)}{6\sqrt{3}\, |F|^2}\,,\nonumber \\
X_8&=&\frac{i}{6\sqrt{3}\,|c|^2}|F|^{-2e^{i\psi}}(\bar{c}^2\bar{F}-
c^2{F}|F|^{2i\sqrt{3}})\,,
\label{rdvecforcompfornonsplitforexamp}
\end{eqnarray}
where $\psi$ is given in (\ref{equationforpsi}) and $c$ is a 
complex constant. The corresponding first fundamental form is 
immediately obtained from (\ref{indmetricfornonsplitting}) for 
$N=3$ and given as
\begin{equation}
I=\frac{2}{3}\, \frac{|F^{\prime}|^2}{|F|^{2}}\,d\xi\,d\bar{\xi}\,.
\label{firstfunforcp2fornonsplitexamp}
\end{equation}
Note that the components of the radius vector $\vec{X}$ in 
(\ref{rdvecforcompfornonsplitforexamp}) satisfy the following relations
\begin{equation}
X_1^2+X_2^2=X_5^2+X_7^2=X_6^2+X_8^2=\frac{1}{27}\,. 
\end{equation}
Eliminating the functions $F$ and $\bar{F}$ in 
(\ref{rdvecforcompfornonsplitforexamp}) we obtain
\begin{eqnarray}
X_1 &=& \frac{i}{6\sqrt{3}\,|c|^2}e^{-(v+\bar{v}) e^{i\psi}}(\bar{c}^2 e^v-c^2 e^{\bar{v}} e^{i \sqrt{3} (v+\bar{v})})\,, \nonumber \\
X_2 &=& -\frac{1}{6\sqrt{3}\,|c|^2}e^{-(v+\bar{v}) e^{i\psi}}(\bar{c}^2 e^v+c^2 e^{\bar{v}} e^{i \sqrt{3} (v+\bar{v})})\,, \nonumber \\
X_5 &=& -\frac{1}{3\sqrt{3}}\cos\left(\frac{3}{2}(\sqrt{3}X_3+X_4)\right)\,, \nonumber \\
X_6 &=& \frac{1}{6\sqrt{3}\,|c|^2}e^{-(v+\bar{v}) e^{i\psi}}(\bar{c}^2 e^{\bar{v}}+c^2 e^{{v}} e^{i \sqrt{3} (v+\bar{v})})\,, \nonumber \\
X_7 &=& -\frac{1}{3\sqrt{3}}\sin\left(\frac{3}{2}(\sqrt{3}X_3+X_4)\right)\,, \nonumber \\
X_8 &=& \frac{i}{6\sqrt{3}\,|c|^2}e^{-(v+\bar{v}) e^{i\psi}}(\bar{c}^2 e^{\bar{v}}-c^2 e^{{v}} e^{i \sqrt{3} (v+\bar{v})})\,,
\end{eqnarray}
where $\psi$ is given in (\ref{equationforpsi}) and $v=\frac{3}{4} (1+i\sqrt{3}) (X_3+iX_4)$. The surface is parametrized in terms of $X_3$ and $X_4$. Now, the corresponding first fundamental form becomes
\begin{equation}
I=\frac{3}{2}(dX_3^2+dX_4^2)\,.
\end{equation}
Note that this is just the real form of (\ref{firstfunforcp2fornonsplitexamp}) when $\xi^1=X_3$ and $\xi^2=X_4$.

\section{Summary and concluding remarks}
The objective of this paper was to revise and expand on theoretical 
results in \cite{Grundlandetal} concerning surfaces related to the 
$\field{C}P^{N-1}$ sigma model.
In addition, we gave a comprehensive summary of geometric properties and 
corrected mistakes in \cite{Grundlandetal}. 
For example, Proposition $4$ in \cite{Grundlandetal} concerning the 
structural equations for the $\field{C}P^2$ model (where only the 
holomorphic solutions were assumed), has been restated as Proposition $2$.
In doing so, we covered in greater detail the geometrical aspects of 
surfaces immersed in the $su(N)$ algebra. 
Furthermore, we have derived the formulae in terms of explicit functions 
in the $\field{C}P^{N-1}$ model, which makes the results in 
\cite{Grundlandetal} more transparent and useful.

We also computed the Lie-point symmetries of the $\field{C}P^{N-1}$ model 
equations for arbitrary $N$.  
The resulting symmetry algebra is decomposed as a direct sum of two 
infinite-dimensional simple Lie algebras and the $su(N)$ algebra. 
Using the Lie-point symmetries, the method of symmetry reduction can 
now be applied to find solutions which are invariant under subgroups of 
$SU(N)$ with generic orbits of codimension one. 
In \cite{Grundland3}, this analysis was carried out for $N=2$. 
The obtained invariant solutions are complicated expressions in 
terms of elliptic functions. 
As was shown in \cite{Grundland3}, for some cases the reduced ordinary 
differential equations (ODEs) can be transformed into 
the standard form of the P3 Painlev\'{e} transcendent. 
Matters get worse when  $N\geq 3.$
Although the reduction can still be carried out, the resulting ODEs are 
coupled and do not appear to be separable. 
One can prove the existence of solutions but `how to find them' remains an 
open problem.

For the $\field{C}P^2$ model, we characterized the immersion of surfaces 
in the $su(3)$ algebra. 
Explicit formulae were found for the moving frame, the structural equations 
(Gauss-Weingarten and Gauss-Codazzi), the first and second fundamental forms, 
the Gaussian, the mean curvatures, the Willmore functional 
and the topological charge. 
These quantities are expressed in terms of holomorphic fields of the 
$\field{C}P^2$ model. 
The theoretical concepts have been illustrated with various examples. 
We also have shown that non-degenerate affine surfaces in $\field{R}^8$ 
associated with the $\field{C}P^2$ model are affine spheres. 
Finally, we discussed dilation-invariant solutions of the 
$\field{C}P^{N-1}$ model, holomorphic immersion of surfaces associated 
with $\field{C}P^2$ models, and mixed soliton-type solutions of the 
$\field{C}P^2$ model and its corresponding surfaces.

\section*{ACKNOWLEDGMENTS}
This work is supported in part by research grants from NSERC of Canada. 
W.\ H.\ gratefully acknowledges the financial support and hospitality of 
the CRM during his sabbatical leave. 
\.{I}.Y. acknowledges a postdoctoral fellowship 
awarded by the Laboratory of Mathematical Physics of the CRM, 
Universit\'{e} de Montr\'{e}al.

\vfill
\newpage
\section*{Appendix A \label{normals}}
In this Appendix we give the explicit form of the vector normals, 
\begin{eqnarray}
\eta_j=\phi^{\dagger} s_j\phi\,, \qquad j=3,\ldots, 8\,,
\nonumber
\end{eqnarray}
to the surface immersed in the $su(3)$ algebra. 
The general expressions are too complicated to be useful.
Instead, we consider the case of a $2D$ surface associated with the 
$\field{C}P^2$ model with solution (\ref{holsolutionscp2special}).

We present the normals in the equivalent matrix form. 

The first normal is 
\begin{eqnarray}
\eta_3=\phi^{\dagger} s_3\phi = i\eta^3_{ij}\,,
\nonumber
\end{eqnarray}
where
\begin{eqnarray}
&&\eta^3_{11}=
\frac{4(|\xi|^2-1)}{{\Gamma_2}^{\!\!2}}\,,
\,\,
\eta^3_{12}=
\frac{2\xi\big(4+|\xi|^2\Gamma_1\big)}{\Gamma_1{\Gamma_2}^{\!\!2}}\,,
\,\,
\eta^3_{13}=
\frac{2\xi^2\Gamma_5}{\Gamma_1{\Gamma_2}^{\!\!2}}\,, \nonumber \\
&&\eta^3_{21}=
\frac{2\bar{\xi}\big(4+|\xi|^2 \Gamma_1 \big)}{\Gamma_1{\Gamma_2}^{\!\!2}} \,,
\quad 
\eta^3_{22}=
\frac{4+
|\xi|^4 \big(5+|\xi|^2 \Gamma_2 \big)}{{\Gamma_1}^{\!\!2}{\Gamma_2}^{\!\!2}}\,,
\nonumber \\
&&\eta^3_{23}=
-\frac{4{\xi}(|\xi|^2-1)}{{\Gamma_1}^{\!\!2}{\Gamma_2}^{\!\!2}} \,, \quad
\eta^3_{31}=\frac{2\bar{\xi}^2 \Gamma_5}{\Gamma_1{\Gamma_2}^{\!\!2}} \,,
\nonumber \\
&&\eta^3_{32}=
-\frac{4\bar{\xi}(|\xi|^2-1)}{{\Gamma_1}^{\!\!2}{\Gamma_2}^{\!\!2}} \,,
\quad
\eta^3_{33}=
\frac{|\xi|^2\big(4-
|\xi|^2 {\Gamma_3}^{\!\!2}\big)}{{\Gamma_1}^{\!\!2}{\Gamma_2}^{\!\!2}}
\,, \label{norapp1}
\end{eqnarray}
with $\Gamma_j$ ($j=1,\ldots, 5$) defined as
\begin{equation}
\Gamma_j=j+|\xi|^2\,, \qquad j=1,\ldots, 5\,.
\end{equation}
The second normal is
\begin{eqnarray}
\eta_4=\phi^{\dagger} s_4\phi = i\eta^4_{ij}\,,
\nonumber
\end{eqnarray}
where 
\begin{eqnarray}
&&\eta^4_{11}=
\frac{2\big(2+|\xi|^2(2-|\xi|^2)\big)}{\sqrt{3}\,{\Gamma_2}^{\!\!2}}\,,
\quad
\eta^4_{12}=
\frac{2\sqrt{3}\,|{\xi}|^2 {\xi}}{{\Gamma_2}^{\!\!2}}\,,
\nonumber \\
&&\eta^4_{13}=
-\frac{2\sqrt{3}\, {\xi}^2}{{\Gamma_2}^{\!\!2}}\,,
\quad
\eta^4_{21}=
\frac{2\sqrt{3}\, |{\xi}|^2 \bar{\xi}}{{\Gamma_2}^{\!\!2}} \,,
\nonumber \\
&&\eta^4_{22}=
\frac{4+|\xi|^2(|\xi|^2-8)}{\sqrt{3}\,{\Gamma_2}^{\!\!2}}\,,
\quad
\eta^4_{23}=
\frac{4\sqrt{3}\, {\xi}}{{\Gamma_2}^{\!\!2}}\,,
\nonumber \\
&&\eta^4_{31}=
-\frac{2\sqrt{3}\, \bar{\xi}^2}{{\Gamma_2}^{\!\!2}}\,,
\quad
\eta^4_{32}= 
\frac{4\sqrt{3}\, \bar{\xi}}{{\Gamma_2}^{\!\!2}} \,,
\nonumber \\
&&\eta^4_{33}=
\frac{|\xi|^2 \Gamma_4-8}{\sqrt{3}\,{\Gamma_2}^{\!\!2}}
\,. \label{norapp2}
\end{eqnarray}
The next one is
\begin{eqnarray}
\eta_5=\phi^{\dagger} s_5\phi = ie^{-\frac{3i\varphi}{2}}\eta^5_{ij}\,,
\nonumber
\end{eqnarray}
where 
\begin{eqnarray}
&&\eta^5_{11}=
\frac{2|\xi|(e^{3i\varphi}\xi^2-\bar{\xi}^2)}{{\Gamma_2}^{\!\!2}} \,,
\quad
\eta^5_{12}=
-\frac{\sqrt{\xi}\big(4e^{3i\varphi} {\xi}^2\Gamma_1 + \bar{\xi}^2 
(2+|\xi|^2 \Gamma_1)\big)}{\sqrt{\bar{\xi}}\, \Gamma_1{\Gamma_2}^{\!\!2}} \,,
\nonumber \\
&&\eta^5_{13}=
\frac{2{\xi}^{(3/2)}(2e^{3i\varphi} {\xi}\Gamma_1-
\bar{\xi}^3)}{{\bar{\xi}^{(3/2)}}\,\Gamma_1{\Gamma_2}^{\!\!2}}\,,
\quad
\eta^5_{21}=
\frac{\sqrt{\bar{\xi}}\big(4 \bar{\xi}^2\Gamma_1+e^{3i\varphi} {\xi}^2
(2+|\xi|^2 \Gamma_1 \big)}{\sqrt{{\xi}}\, \Gamma_1{\Gamma_2}^{\!\!2}}\,,
\nonumber \\
&&\eta^5_{22}=
-\frac{2(e^{3i\varphi}\xi^2-\bar{\xi}^2)
(2+|\xi|^2 \Gamma_1)}{|\xi|\Gamma_1{\Gamma_2}^{\!\!2}} \,,
\quad
\eta^5_{23}=
\frac{2\sqrt{{\xi}}\big(2\bar{\xi}^3+ e^{3i\varphi}\xi 
(2+|\xi|^2 \Gamma_1)\big)}{{\bar{\xi}^{(3/2)}}\Gamma_1{\Gamma_2}^{\!\!2}}\,,
\nonumber \\
&&\eta^5_{31}=
\frac{2\bar{\xi}^{(3/2)}\big(e^{3i\varphi}\xi^3-
2 \bar{\xi}\Gamma_1 \big)}{{{\xi}^{(3/2)}}\Gamma_1{\Gamma_2}^{\!\!2}}\,,
\quad
\eta^5_{32}=
-\frac{2\sqrt{\bar{\xi}}\big(2e^{3i\varphi}\xi^3+\bar{\xi}
(2+|\xi|^2 \Gamma_1)\big)}{{{\xi}^{(3/2)}}\Gamma_1{\Gamma_2}^{\!\!2}}\,,
\nonumber \\
&&\eta^5_{33}=\frac{4(e^{3i\varphi}\xi^2-
\bar{\xi}^2)}{|\xi|\Gamma_1{\Gamma_2}^{\!\!2}}
\,. \label{norapp3}
\end{eqnarray}
Normal $\eta_6$ is given by
\begin{eqnarray}
\eta_6=\phi^{\dagger} s_6\phi = ie^{-\frac{3i\varphi}{2}}\eta^6_{ij}\,,
\nonumber
\end{eqnarray}
where 
\begin{eqnarray}
&&\eta^6_{11}=
-\frac{2|\xi|\big(e^{3i\varphi}\xi-\bar{\xi}\big)}{{\Gamma_2}^{\!\!2}}\,,
\quad
\eta^6_{12}=
\frac{2{\xi}^{(3/2)}\big(2e^{3i\varphi}\Gamma_1-
\bar{\xi}^2\big)}{\sqrt{\bar{\xi}}\,\Gamma_1{\Gamma_2}^{\!\!2}}\,,
\nonumber \\
&&\eta^6_{13}=
-\frac{{\xi}^{(3/2)}\big(4e^{3i\varphi}\Gamma_1+|{\xi}|^2
\bar{\xi}^2\Gamma_3\big)}{{\bar{\xi}^{(3/2)}}\Gamma_1{\Gamma_2}^{\!\!2}}\,,
\quad
\eta^6_{21}=
-\frac{2{\bar{\xi}^{(3/2)}}\big(2-e^{3i\varphi}{\xi}^2+
2|\xi|^2\big)}{\sqrt{{\xi}}\,\Gamma_1{\Gamma_2}^{\!\!2}}\,,
\nonumber \\
&&\eta^6_{22}=
-\frac{4|\xi|\big(e^{3i\varphi}\xi-
\bar{\xi}\big)}{\Gamma_1{\Gamma_2}^{\!\!2}}\,,
\quad
\eta^6_{23}=
\frac{2{\xi}^{(3/2)}\big(2 e^{3i\varphi} + 
\bar{\xi}^2 \Gamma_3\big)}{\sqrt{\bar{\xi}}\,\Gamma_1{\Gamma_2}^{\!\!2}}\,,
\nonumber \\
&&\eta^6_{31}=
\frac{{\bar{\xi}^{(3/2)}}\big(4+4|\xi|^2+e^{3i\varphi} |{\xi}|^{2} 
{\xi}^2 \Gamma_3\big)}{{{\xi}^{(3/2)}}\Gamma_1{\Gamma_2}^{\!\!2}}\,,
\quad
\eta^6_{32}=
-\frac{2{\bar{\xi}^{(3/2)}}\big(2+e^{3i\varphi} {\xi}^{2} 
\Gamma_3 \big)}{\sqrt{{\xi}}\,\Gamma_1{\Gamma_2}^{\!\!2}}\,,
\nonumber \\
&&\eta^6_{33}=
\frac{2|\xi|\big(e^{3i\varphi}\xi-\bar{\xi}\big)
\Gamma_3}{\Gamma_1{\Gamma_2}^{\!\!2}}
\,. \label{norapp4}
\end{eqnarray}
Normal $\eta_7$ is given by 
\begin{eqnarray}
\eta_7=\phi^{\dagger} s_7\phi = e^{-\frac{3i\varphi}{2}}\eta^7_{ij}\,,
\nonumber
\end{eqnarray}
where
\begin{eqnarray}
&&\eta^7_{11}=
-\frac{2|\xi|\big(e^{3i\varphi}\xi^2+\bar{\xi}^2\big)}{{\Gamma_2}^{\!\!2}}\,,
\quad
\eta^7_{12}=
\frac{\sqrt{{\xi}}\big(4 e^{3i\varphi} {\xi}^{2} \Gamma_1-\bar{\xi}^{2}
(2+|\xi|^2\Gamma_1)\big)}{\sqrt{\bar{\xi}}\,\Gamma_1{\Gamma_2}^{\!\!2}}\,,
\nonumber \\
&&\eta^7_{13}=
-\frac{2{\xi}^{(3/2)} \big(\bar{\xi}^3+2 e^{3i\varphi} 
\xi\Gamma_1 \big)}{{\bar{\xi}^{(3/2)}}\Gamma_1{\Gamma_2}^{\!\!2}}\,,
\quad
\eta^7_{21}=
\frac{\sqrt{\bar{\xi}} \big(4 \bar{\xi}^2 \Gamma_1 - e^{3i\varphi} 
{\xi}^2 (2+|\xi|^2\Gamma_1)\big)}{\sqrt{{\xi}}\,\Gamma_1{\Gamma_2}^{\!\!2}}\,,
\nonumber \\
&&\eta^7_{22}=
\frac{2\big(e^{3i\varphi}\xi^2+\bar{\xi}^2\big) 
\big(2+|\xi|^2\Gamma_1\big)}{|\xi|\Gamma_1{\Gamma_2}^{\!\!2}}\,,
\quad
\eta^7_{23}=
\frac{2\sqrt{\xi} \big(2 \bar{\xi}^3 - e^{3i\varphi}\xi 
(2+|\xi|^2\Gamma_1) \big)}{{\bar{\xi}^{(3/2)}}\Gamma_1{\Gamma_2}^{\!\!2}}\,,
\nonumber \\
&&\eta^7_{31}=
-\frac{2\bar{\xi}^{(3/2)} \big(e^{3i\varphi}\xi^3+
2 \bar{\xi} \Gamma_1 \big) }{{{\xi}^{(3/2)}}\Gamma_1{\Gamma_2}^{\!\!2}}\,,
\quad
\eta^7_{32}=
\frac{2 \sqrt{\bar{\xi}} \big(2e^{3i\varphi}\xi^3-\bar{\xi}
(2+|\xi|^2\Gamma_1) \big)}{{{\xi}^{(3/2)}}\Gamma_1{\Gamma_2}^{\!\!2}}\,,
\nonumber \\
&&\eta^7_{33}=
-\frac{4\big(e^{3i\varphi}\xi^2+
\bar{\xi}^2\big)}{|\xi|\Gamma_1{\Gamma_2}^{\!\!2}}
\,. \label{norapp5}
\end{eqnarray}
The last normal is given by 
\begin{eqnarray}
\eta_8=\phi^{\dagger} s_8\phi = e^{-\frac{3i\varphi}{2}}\eta^8_{ij}\,,
\nonumber
\end{eqnarray}
where
\begin{eqnarray}
&&\eta^8_{11}=
\frac{2|\xi|\big(e^{3i\varphi}\xi+\bar{\xi}\big)}{{\Gamma_2}^{\!\!2}}\,,
\quad
\eta^8_{12}=
-\frac{2{\xi}^{(3/2)}\big(\bar{\xi}^2+
2e^{3i\varphi}\Gamma_1\big)}{\sqrt{\bar{\xi}}\,\Gamma_1{\Gamma_2}^{\!\!2}}\,,
\nonumber \\
&&\eta^8_{13}=
\frac{{\xi}^{(3/2)}\big(4e^{3i\varphi}\Gamma_1-|{\xi}|^2\bar{\xi}^2
\Gamma_3\big)}{{\bar{\xi}^{(3/2)}}\Gamma_1{\Gamma_2}^{\!\!2}}\,,
\quad
\eta^8_{21}=
-\frac{2{\bar{\xi}^{(3/2)}}\big(2+e^{3i\varphi}{\xi}^2+
2|\xi|^2\big)}{\sqrt{{\xi}}\,\Gamma_1{\Gamma_2}^{\!\!2}}\,,
\nonumber \\
&&\eta^8_{22}=\frac{4|\xi|\big(e^{3i\varphi}\xi+
\bar{\xi}\big)}{\Gamma_1{\Gamma_2}^{\!\!2}}\,,
\quad
\eta^8_{23}=-\frac{2{\xi}^{(3/2)}\big(2 e^{3i\varphi} - 
\bar{\xi}^2 \Gamma_3\big)}{\sqrt{\bar{\xi}}\,\Gamma_1{\Gamma_2}^{\!\!2}}\,,
\nonumber \\
&&\eta^8_{31}=\frac{{\bar{\xi}^{(3/2)}}\big(4+4|\xi|^2-
e^{3i\varphi} |{\xi}|^{2} {\xi}^2 \Gamma_3\big)}{{{\xi}^{(3/2)}}
\Gamma_1{\Gamma_2}^{\!\!2}}\,,
\quad
\eta^8_{32}=
\frac{2{\bar{\xi}^{(3/2)}}\big(e^{3i\varphi} {\xi}^{2} 
\Gamma_3 -2\big)}{\sqrt{{\xi}}\,\Gamma_1{\Gamma_2}^{\!\!2}}\,,
\nonumber \\
&&\eta^8_{33}=
-\frac{2|\xi|\big(e^{3i\varphi}\xi+
\bar{\xi}\big)\Gamma_3}{\Gamma_1{\Gamma_2}^{\!\!2}}
\,. \label{norapp6}
\end{eqnarray}

\end{document}